\documentclass[a4paper,12pt,reqno]{amsart}
\usepackage{amssymb}\theoremstyle{plain}

\usepackage{amscd}
\usepackage{enumitem}

\newtheorem{theorem}{Theorem}
\newtheorem*{theorem*}{Theorem}
\newtheorem{corollary}{Corollary}
\newtheorem*{corollary*}{Corollary}
\newtheorem{lemma}{Lemma}
\newtheorem*{lemma*}{Lemma}
\newtheorem{proposition}{Proposition}
\newtheorem*{proposition*}{Proposition}
\newtheorem{conjecture}{Conjecture}
\newtheorem*{conjecture*}{Conjecture}
\theoremstyle{definition}
\newtheorem{definition}{Definition}
\newtheorem*{definition*}{Definition}
\theoremstyle{remark}
\newtheorem{remark}{Remark}
\newtheorem*{remark*}{Remark}
\newtheorem{example}{Example}
\newtheorem*{problem*}{Problem}
\newtheorem{problem}{Problem}

\def\mathbi#1{\textbf{\em #1}}

\begin{document}

\begin{center}
\title[ArAlg]{On the arithmetic and algebraic properties of Minkowski balls and spheres}
\end{center}

\maketitle

\begin{center}
{\bf Nikolaj M. Glazunov } \end{center}

\begin{center}
{\rm Glushkov Institute of Cybernetics NASU, Kiev, } \\

{\rm  Email:} {\it glanm@yahoo.com }
\end{center} 

\bigskip


\bigskip


\footnote{The author was partially supported by Simons grant 992227.}  

\begin{abstract}
This paper gives a brief overview of some new work in number theory
 and algebra, and also studies the arithmetic and algebraic properties 
of Minkowski balls and spheres. The content of the paper is presented in more detail in the table of Contents and in the Introduction.
\end{abstract}

\newpage

\tableofcontents

\newpage

\section{Introduction}

Problems of studying lattice points on spheres and on the boundaries of 
centrally symmetric  convex bodies of variable radius belong to classical
 problems of number theory.
 The results of the study of such problems, in addition to number theory, 
find applications in coding theory, in combinatorics, and in other applications.
To understand varieties and manifolds,   varieties and manifolds with additional structures, 
as well as as well as problems associated with these objects,  it is necessary to understand the symmetries and invariants of these objects.
 The two-dimensional Minkowski ball is associated with its critical lattice, and its invariant is its critical determinant (see section \ref{Minkowski balls}).
 For a manifold in $n$-dimensional Euclidean space ${\mathbb R}^n$,
 one studies  the group action on the manifold by means of diffeomorphisms.
 An example of such an action is the action of the finite subgroup of orthogonal group ${\mathbi O}^n$
  on the unit sphere ${\mathbb S}^{n-1}$.
  If the manifold has a cellular decomposition, the group is finite and the 
action is smooth, then the action of the group is given by permutations of cells.
For instance if the action of the group $G$ is free and ${\mathbb Z}G$ 
is the integral group ring, then the cellular chain complex is the chain complex of free ${\mathbb Z}G$- modules.

In linear algebra, in field theory, in topology, in probability theory, in graph theory, and also in a number of applied areas, for example, in coding theory, there are concepts of independence.

Along with traditional domains (fields and classical rings) of definition of manifolds and schemes, absolute objects are currently being actively studied as domains of definition.
A. Connes and C. Consani, as well as other researchers
 (see \cite{ccaa,ccag,ccsp,ccrrz,ccars,ccaar} and references therein),
 define and study absolute algebra,  absolute algebraic
 geometry (affine case) and approaches to absolute arithmetic.
 
 It is possible that the transition from the above-mentioned 
traditional areas to the absolute theories of these is a 
transition to the combinatorial structure of these traditional 
areas.

Since Minkowski spheres can be both algebraic and transcendental, it is interesting to explore their dynamics and compare them. We will begin by examining the elements of algebraic dynamics in this case.
 
 Whitney's matroid theory \cite{whi} and its development are sometimes 
referred to as independence theory.

In view of this, and also to study the combinatorial structures of the objects under study, we present a (very unsystematic) review of matroid theory.

This paper gives a brief overview of some new works in these directions, 
and also studies the arithmetic and algebraic properties of Minkowski balls 
and spheres.

The plan of the paper is as follows.

\begin{itemize}
  
\item  {\it  Introduction}

\item {\it Minkowski balls and spheres}

\item {\it About the average  values of number-theoretic objects}

\item {\it Packing of  Minkowski balls}

\item {\it On  coverings by Minkowski balls}

\item {\it Critical lattices, Minkowski domains and  their optimal packing}

\item {\it Inscribed and circumscribed hexagons of minimum areas}

\item {\it  On coverings by Minkowski domains in the plane}

\item {\it  Direct systems,  direct limits and duality}






\item {\it About Riemann-Roch theorem for algebraic curves,  lattices and arithmetic curves}

\item {\it Algebraic dynamics}

\item {\it  On elements of matroid theory}

\item {\it Minkowski  spheres and matroids
}

\end{itemize}



\section{Minkowski balls and spheres}

  Let $x_1, \cdots x_n  \in {\mathbb R}$ and let 
 \begin{equation}
\label{smt}   D_p: \;  |x_1|^p + \cdots + |x_n|^p < 1, \; p \ge 1
\end{equation}
be ($n$-dimensional) Minkowski's balls with boundary 
\begin{equation}
\label{mt}
{{\mathbb S}^{n-1}_c}: \;   |x_1|^c + \cdots + |x_n|^c  = 1, \;  c \ge 1.
\end{equation}
 (Minkowski's spheres).\\
\begin{remark}
Keeping the previously introduced notation $D_p$ for Minkowski balls, for n-dimensional Minkowski spheres we use the notation ${{\mathbb S}^{n-1}_c}$ where $ c \ge 1$.
The unit Euclidean $n$-dimensional sphere denotes as
 ${\mathbb S}^{n-1} = {{\mathbb S}^{n-1}_2}$. 
\end{remark}

\begin{remark}
The cases $p=1$ and $p \to \infty$ are called the limit $n$-dimensional Minkowski balls.
Accordingly, for $c = 1$ and $c \to \infty$ we have limiting $n$-dimensional Minkowski spherses.
\end{remark}
Let
\begin{equation}
\label{ams}
 \{\lfloor|x_1|^c\rfloor + \{\lfloor|x_2|^c\rfloor + \cdots +   \{\lfloor|x_n|^c\rfloor = m
  \}
\end{equation} 
 be the arithmetic  Minkowski's spheres with any $m \in {\mathbb Z}_+$ and $c > 1$. 
\begin{remark*}
 Equations similar to (\ref{ams}) were also considered by Segal \cite{seg} when studying deformations of Waring’s problem. See also
 \cite{de,az,ilms}. 
\end{remark*}
 We will also consider complex Minkowski's balls of the form
 \begin{equation}
\label{cmt}
 |x_1|^p + \cdots + |x_s|^p  + |y_1 + z_1|^{p/2} + \cdots + |y_t + z_t|^{p/2} < 1, \; p \ge 1
 \end{equation}
 
Arithmetic also considers  algebraic extensions $K/{\mathbb Q}$ of degree $n$ of the field ${\mathbb Q}$. Minkowski interpret points of $K$ as points of $n$-dimensional space \cite{Mi:DA,neu}.This leads to Minkowski space in algebraic number theory.

\section{About the average  values of number-theoretic objects 
  }
  Many number-theoretic objects, for example, arithmetic functions, 
trigonometric sums, and others, fluctuate significantly over the domain of
 definition.
 The well-known method of I.M. Vinogradov's \cite{vi} estimates of H. Weil's trigonometric sums were successfully applied to solve fundamental problems of \cite{vi1,vi2} number theory.
A.A. Karatsuba generalized the method of I.M. Vinogradov, reducing by smoothing the estimate of the individual trigonometric sum of H. Weil to an estimate of its average value \cite{akc}.
While studying the asymptotically uniform distribution of integer points on  three-dimensional spheres of increasing radii, 
Linnik \cite{lin} discovered the connection between this question with the ergodic theory of abstract discrete-time flow.
 Ergodic, discrete and  lacunary (spherical) averages have been investigated in papers by  Anderson, Cook, Hughes, Kumchev  \cite{achk}, by  Bourgain, Mirek, Stein,  Wróbel \cite{bmsw}, by Cook, Hughes \cite{ch}, by  Iosevich,  Sawyer,  Seeger
  \cite{iss}, by Magyar, Stein,  Wainger \cite{msw}.
  
  Let ${{\mathbf S}^{3}_c} := \{ x \in {\mathbb Z}^3 : \{\lfloor|x_1|^c\rfloor + \{\lfloor|x_2|^c\rfloor +  \{\lfloor|x_3|^c\rfloor =\lambda \}$  be the arithmetic sphere with any $\lambda \in {\mathbb Z}_{+}$ and $c > 1$.
  In the paper \cite{ilms}  ergodic averages $A^{c,3}_{\lambda} f(x)$ over arithmetic spheres and discrete averages 
$M^{c,3}_{\lambda} f(x)$ on $({\mathbb Z}^3, \nu_{\mathbb Z}^3, T_{\mathbb Z}^3)$ are considered, where $c \in (1, 11/10)$, and in the first case for  $\sigma$-finite measure space $(X, \nu), f \in L^0 (X)$,  and in the second case  
$f : {\mathbb Z}^3 \to \hat{\mathbb C} $ .
Authors  investigate the lattice points on the arithmetic sphere 
$ {\mathbf S}_{h_1,h_2,h_3} (\lambda) := \{ x \in ({\mathbb Z} \setminus [N_0 - 1])^3 : \lfloor|h_1(x)| \rfloor + \lfloor|h_2(x)|\rfloor +  \lfloor|h_3 (x)| \rfloor =\lambda \}$ where  $N_0 \in {\mathbb Z}_{+}$ is a sufficiently large absolute constant and functions $h_1, h_2, h_3$ are constant multiples of regularly varying functions. 
``Continuous spherical averaging operators are defined by'' 
 ${\mathcal A}_t^c f(x) := \int_{{\mathbb S}^2_c} f(x - t\theta) d{\mu}_c(\theta), x \in {\mathbb R}^3, t > 0, 
 f \in {\mathcal C}^{\infty}_c ({\mathbb R}^3),  $
 where ${\mu}_c$ is the measure on ${\mathbb S}^2_c$ and ${\mathcal C}^{\infty}_c ({\mathbb R}^3)$ is the set of all compactly supported smooth functions on ${\mathbb R}^3.$
  By analogy with  equidistribution problem by Linnik \cite{lin}  the authors study projections
 ${{\mathbf P}^{3}_c}(\lambda) := \{\lambda^{-1/c} x: x \in  {\mathbf S}^{3}_c (\lambda) \}$  in the neighborhood of the unit (Minkowski) sphere ${\mathbb S}^{2}_c, c \in (1, 9/8)$.
 Let  $r_c (\lambda) = \#{\mathbf S}^{3}_c$ be the number of lattice points in ${\mathbf S}^{3}_c$ and
  $\phi \in {\mathcal C}^{\infty} ({\mathbb R}^3)$.
The equidistribution problem is reduced to proving 
  $\frac{1}{r_c (\lambda)} \sum_{x \in {{\mathbf P}^{3}_c}(\lambda)} \phi(x) \to \int_{{\mathbb S}^{2}_c} \phi(x) d\nu_c(x)$ as
 $\lambda \to \infty$.
 Here $\nu_c = \frac{\mu_c}{\mu_c({\mathbb S}^{2}_c)} $ 
is a probability measure on 
${{\mathbb S}^{2}_c}$. 
 The discrepancy with respect to caps measures the irregularity of distribution of lattice points on spheres.
The paper proves the following theorems.
The main result is Theorem 1.7: The maximal inequalities for full discrete averages and for lacunary discrete averages  of $M^{c,3}_{\lambda} f$.

 Theorem 1.5: The ergodic theorem for averages.
 Theorem 1.5  is the second main result of the paper.
Theorem 1.2: Asymptotic formula for number of  lattice points on ${\mathbf S}_{h_1,h_2,h_3} (\lambda)$.
Theorem 1.8: Maximal and $r$-variational estimates for the operators ${\mathcal A}_t^c$.
Theorem 1.9: Equidistribution of the points from ${{\mathbf P}^{3}_c}(\lambda)$ on the unit sphere ${\mathbb S}^{2}_c$.
Theorem 1.10: Estimation of discrepancy.
Let us briefly sketch the proof of the main results.
In proving  Theorem 1,7 the authors make use ideas from the circle method and for every $\lambda \in {\mathbb Z}_{+}$ split the interval $[-1/2, 1/2]$ into major and minor arcs.
The authors then carry out estimates on minor and  major arcs.
The maximal functions corresponding to the major arcs are estimated.
Values ${\mu_c({\mathbb S}^{2}_c)} = \frac{8\Gamma(1/c)^3}{c^2\Gamma(3/c)}$ are given.
Estimates of the Fourier transform for measures $\mu_c$ associated with spheres ${\mathbb S}^{2}_c$ are the longest and most technical part of the paper.
The maximum ergodic part of Theorem 1.5 follows from the maximum inequality  for full averages from Theorem 1.7 by Calderon transference principle as discribed by Bourgain \cite{bou}.
Accordingly, the lacunary variant of the inequality implies the maximal inequality for the lacunary case of Theorem 1.5 by the above versus of Calderon transference  principle.

\section{Packing of  Minkowski balls} 
\label{Minkowski balls}
 We investigate lattice packings of  Minkowski balls. By the results of the proof of Minkowski conjecture about the critical determinant 
we devide Minkowski balls on 3 classes: Minkowski balls, Davis balls and  Chebyshev-Mordell balls \cite{gl1}.
\begin{remark}
 Here and below we use the phrase Chebyshev-Mordell (instead 
of the previously used phrase Chebyshev-Cohn). The phrase 
Chebyshev-Mordell better corresponds to the history of the 
study of the Minkowski hypothesis.
\end{remark}

We investigate lattice packings of  these balls on  planes  and  search among these packings the optimal packings. In the paper \cite{gl1} we prove that the optimal lattice packing of the Minkowski, Davis, and Chebyshev-Mordell balls is realized with respect to the sublattices of index two of the critical lattices of
corresponding  balls.

A system of equal  balls in $n$-dimensional space is said to form a packing, if no two balls of the 
system have any inner points in common.

Lately the remarkable results in resolving the problem of optimal packing of Euclidean balls in $8$ and 
$24$-dimensional real Euclidean spaces have been obtained \cite{via,ckmrv}.

\subsection{Minkowski balls and the density  of the packing of $2$-dimensional Minkowski  balls}

\begin{definition}
\label{dc}
Let $p_{0} \in {\mathbb R}$ be the Davis constant such that $2,57 < p_{0}  < 2,58$. \\
\end{definition}

We consider balls of the form
\[
   D_p: \;  |x|^p + |y|^p  < 1 \; p \ge 1,
\]
and call balls with the conditions
\[
|x|^p + |y|^p < 1,\;  2 > p > 1, 
\]
 the  Minkowski balls in two dimension and 
correspondingly the circles with  the condition
\[ 
|x|^p + |y|^p = 1, \; 2 > p > 1 
\]
the  Minkowski circles in two dimension.   
  
  Limiting Minkowski circle  in two dimension:    $|x| + |y| = 1$   \\

Davis balls  in two dimension: 
        $|x|^p + |y|^p <  1$  for $p_{0} > p \ge 2$ \\

Davis circles in two dimension: 
    $|x|^p + |y|^p = 1$  for $p_{0} > p \ge 2$ \\

Chebyshev-Mordell balls in two dimension:   $|x|^p + |y|^p < 1 $  for $ p \ge p_{0}$ \\

Chebyshev-Mordell circles in two dimension: $|x|^p + |y|^p = 1$  for $ p \ge p_{0}$ \\

 Limiting Chebyshev-(Mordell) ball in two dimension: 
\[
||x, y||_{\infty} < \max(|x|,|y|).
 \]
 
Recall the definition of a packing lattice \cite{via,ckmrv,Cassels,lek}.
We will give it for $n$-dimensional Minkowski balls $D_p^n$ in ${\mathbb R}^n$.

\begin{definition}
 Let $\Lambda$ be a full lattice in ${\mathbb R}^n$ and $a \in {\mathbb R}^n$. 
 In the case if it is occurs that no two  of balls $\{D_p^n + b, b \in \Lambda + a \}$ have inner points in common, the collection of balls 
   $\{D_p^n + b, b \in \Lambda + a \}$ is called a $(D_p^n,  \Lambda)$-packing, and $\Lambda$ is called a packing lattice of $D_p^n$.
\end{definition}

Recall also that if $\alpha \in {\mathbb R}$ and $D_p^n$ is a ball than $\alpha D_p^n$ is the set of points $\alpha x, x \in D_p^n$.
 

From the considerations of Minkowski and other authors \cite{Mi:DA,Mi:GZ,Cassels,lek,Gl4}, the following statements can be deduced 

Denote by $V(D_p)$ the volume (area) of $D_p$.   $V(D_p) = 4 \frac{(\Gamma(1 + \frac{1}{p}))^2}{\Gamma(1 + \frac{2}{p})}$.

\begin{proposition} \cite{Cassels,lek}.
\label{p1}
A lattice $\Lambda$ is a packing lattice of $D_p$ if and only if  it is admissible lattice for $2 D_p$.
\end{proposition}

\begin{proposition} \cite{Cassels,lek}.
\label{p3}
 The dencity of a $(D_p, \Lambda)$-packing is equal to $V(D_p)/d(\Lambda)$ and it is maximal if $\Lambda$ is critical for  $2 D_p$.
\end{proposition}

  \begin{theorem} \cite{gl1}
  \label{tolp}
The optimal lattice packing of the Minkowski, Davis, and Mordell balls is realized with respect to the sublattices of index two of the critical lattices
  \[
  (1,0)\in\Lambda_{p}^{(0)},\; (-2^{-1/p},2^{-1/p}) \in \Lambda_{p}^{(1)}.
  \]
    \end{theorem}
    
  \section{On  coverings by Minkowski balls}
   
   Lattice coverings in the real plane by Minkowski balls  with inscribed convex symmetric hexagons generated by admissible lattices of Minkowski balls are studied. The covering problem is dual to the packing problem. In different covering problems, this duality, apparently, can manifest itself in different ways. In this paper, duality is understood as the duality between admissible lattices having 3 pairs of points on the boundary $D_p$, and symmetric convex hexagons (possibly a quadrangle) corresponding to these lattices and determining the values of the determinants of the corresponding covering lattices.
   
   We use admissible lattices that have 3 pairs of points on the boundary of $D_p$.
We construct an explicit moduli space (\ref{Ap})  of areas of symmetric convex hexagons (possibly a quadrangle) which determine values of the determinants of corresponding  covering lattices .

This construction and its consequences make it possible to obtain the best currently known lower bounds for covering constants of Minkowski balls. An estimate from above is also given.
These give the corresponding estimation of the covering densities.
The best known value of covering density of Minkowski ball is obtained (Proposition \ref{emcd}). 

\subsection{Two-dimensional Minkowski balls  and hexagons.}
 \subsubsection{ Hexagons and coverings.} 
Every admissible lattice of $D_p$ containing 3 pairs of points on the boundary of $D_p$ defines a hexagon inscribed in $D_p$.
We denote such a hexagon as $ {\mathcal{H}} _p $, and call it a hexagon of the admissible lattice, briefly: {\it al-hexagon}.

 \begin{remark}
In the limiting cases $p=1$ and $p=\infty$ corresponding hexagons (following  Fejer Toth, see \cite{lek} ) $ {\mathcal{H}} _p $ are quadrangles.
 \end{remark}

 By specifying results from \cite{lek} to the case of Minkowski balles and  al-hexagons $ {\mathcal{H}} _p $ we have the next
 \begin{proposition}
 \label{ccdmb}
 Let $\Gamma(D_p)$ be the covering constant of $D_p$.
 Each Minkowski ball $D_p$ contains an al-hexagon  (possibly a quadrangle) $ {\mathcal{H}} _p $ of maximum area.
 If this area is denoted by $\gamma_h(D_p)$ (or for symmetric convex inscribed al-hexagons  
 $ {\mathcal{H}} _p $ by $ a( {\mathcal{H}} _p)$)  then 
  \begin{equation}
 \label{ccmb} 
       \Gamma(D_p) = \gamma_h(D_p).
       \end{equation}  
       For the density of thinnest lattice covering of $D_p$ we get
       \begin{equation}
 \label{dcmb} 
       \vartheta(D_p) = V(D_p)/\gamma_h(D_p).
 \end{equation}  
\end{proposition} 

 \begin{remark}
The expression (\ref{ccmb}) is the the value of the ${\bf covering \; constant}$ of  Minkowski ball $D_p$, and 
 the expression (\ref{dcmb}) is the ${\bf density}$ of the  ${\bf thinnest \; lattice \;covering}$ within the framework of our considerations.
   \end{remark}
   
  \subsection{Inscribed hexagons of admissible lattices  of Minkowski balls $D_p$ and their moduli space. }
 \begin{theorem}
 \label{msa}
The set of areas of the entire family of al-hexagons (possibly a quadrangle) $ {\mathcal{H}}_p, (1 <  p < \infty)$ inscribed in Minkowski balls $D_p$ is parameterized by the function $A(\sigma, p)$ defining the corresponding moduli space
   $\mathbf{A}: $
 \begin{equation}
 \label{Ap}
 A(\sigma, p) = 3(\tau + \sigma)(1 + \tau^{p})^{-\frac{1}{p}}
  (1 + \sigma^p)^{-\frac{1}{p}}, 
 \end{equation}
 In limiting cases we have the following areas 
 of quadrangles: $ p=1, \; a( {\mathcal{H}} _1) = 2; \; p=\infty, \; a( {\mathcal{H}} _{\infty}) = 4$.
  \end{theorem}
  
   \begin{definition}
    Moduli space (\ref{Ap}) is the set of areas of the entire family of al-hexagons (possibly a quadrangle) $ {\mathcal{H}}_p, (1 \le p \le \infty)$ inscribed in Minkowski balls $D_p$. Briefly we will call (\ref{Ap}) the moduli space of hexagons.
  \end{definition}
  
   \begin{remark}
    The function domain of the  function (\ref{Ap}) is
     $$ {\mathcal M}: \; \infty > p > 1, \; 1 \leq \sigma \leq \sigma_{p} =
 (2^p - 1)^{\frac{1}{p}}, $$
     $ {\mathcal M}$  has two limit points $p=1$ and $p=\infty$.
       \end{remark}
  
  \begin{remark}
  The function (\ref{Ap}) is upwardly convex (concave)  in the  real space with Cartesien coordinates     $(\sigma, p, A)$.
  \end{remark}
  
  \subsection{Sections of  moduli space of hexagons.}
        We will solve the problem of the thinnest covering of the real plane by Minkowski balls with hexagons $ {\mathcal{H}} _p$ if we define a curve of maxima on an upwardly convex (concave) surface (\ref{Ap}).
Unfortunately, to date the author knows only a few points of maxima on this curve.   In this regard, we present here a family of surfaces, give  some intersections of (\ref{Ap}) with these surfaces and formulate a (rough) hypothesis about the curve of maxima on (\ref{Ap}).

   \begin{definition}
     Let 
      \begin{equation}
 \label{fsap}  
   \sigma_{\alpha, p} = (2^p - 1)^\frac{1}{\alpha p}, \;  \alpha \ge 1
      \end{equation}  
     be a family of curves parameterized by ${\alpha }$ with $p$ varying from $1$ to $\infty$      
          in 
              real plane  $\{\sigma, p\}$.
   \end{definition}
   
     \begin{remark}
     $\sigma_{\alpha, p} \neq  \sigma_{ p, \alpha}.$ 
    \end{remark}
   
  \begin{remark}
   Let us consider sections orthogonal to the plane $\{\sigma, p\}$
from curves $\sigma_{\alpha,p}$ to the moduli space $\mathbf{A }$. Such sections determine the covering $A\sigma_{\alpha, p}$ of the moduli space $\mathbf{A}$. Thus curves of the family (\ref{fsap})
     define family of curves
          \begin{equation}
 \label{sap}  
     A(\alpha, p): 3(\tau + (2^p - 1)^\frac{1}{\alpha p})(1 + \tau^{p})^{-\frac{1}{p}}
  (1 + ((2^p - 1)^\frac{1}{\alpha p})^p)^{-\frac{1}{p}} 
  \end{equation}  
  on the function range (\ref{Ap}).
   These curves are parametrized by $\alpha$. For each such $\alpha, \; p  $ varies from $ 1$ to   $\infty \;(1 < p < \infty)$.
      \end{remark}
      
\subsection{Estimates and bounds on covering constants. } c
 Here we estimate covering constants of Minkowski balls from above and from below. 
According to the results of Sas and Dovker (see \cite{lek}) in our case of Minkowski balls we have

\begin{proposition}
\label{mhe}
Let  ${\mathcal H}_p$ be the 
 hexagon of maximum area inscribed in Minkowski ball $D_p$.\\
 Then $V({\mathcal H}_p) = \gamma_h(D_p) \ge \frac{3\sqrt{3}}{2\pi}\cdot V(D_p)$
\end{proposition}

\begin{remark}
\label{teamb}
The trivial estimate from above of $V({\mathcal H}_p) $ is given by the area $V(D_p)$ of Minkowski ball:
$$ V({\mathcal H}_p) \le  V(D_p) $$.
\end{remark}

Let estimate covering constants of Minkowski balls $D_p$ from below on the base of results of the proof of Minkowski conjecture.
To do this, we will use the results of the proof of the Minkowski conjecture (now
theorem \cite{GGM:PM}) and the constructed moduli space (\ref{Ap}).

\begin{theorem}
\label{ma}
\begin{equation}
\label{ema}
\min({\bf A}({\mathcal H}_p)) = \left\{
                   \begin{array}{lc}
    3\cdot4^{-1/p} \frac{1 +\tau_p }{1 - \tau_p}, \; 1 \le p \le 2, \; p \ge p_{0},\\
    \frac{3}{2}{\sigma}_{p}, \;  2 \le p \le p_{0};\\
                     \end{array}
                       \right.
                       \end{equation}
\end{theorem}

To estimate the covering constants of Minkowski balls $D_p$ from below, we define the inverse minimum 
$i-\min({\bf A}({\mathcal H}_p)) $ for  ${\bf A}({\mathcal H}_p) $. 

\begin{definition}
\begin{equation}
 \label{ima}
    i-\min({\bf A}) =   i-\min({\bf A}({\mathcal H}_p)) = \left\{
                   \begin{array}{lc}
              \frac{3}{2}{\sigma}_{p}, \;   1 \le p \le 2, \; p \ge p_{0},  ;\\      
    3\cdot4^{-1/p} \frac{1 +\tau_p }{1 - \tau_p}, \; 2 \le p \le p_{0}.\\
                        \end{array}
                       \right.
                       \end{equation}
\end{definition}

\begin{example}
    In the case $p=3$ for moduli space $\mathbf{A}$ (see (\ref{Ap})) we have \\
    $ \min({\bf A}({\mathcal H}_3))     \approx      2.859     $ \\
     $ i-\min({\bf A})({\mathcal H}_3))        \approx     2.870       $\\
     $  \min({\bf A}({\mathcal H}_3)) <  i-\min({\bf A})({\mathcal H}_3)). $
 \end{example}

From theorems \ref{msa}, \ref{ma} with definition  \ref{ma} we have the next estimation from below of the covering constant $\Gamma(D_p)$.

\begin{proposition}
\begin{equation}
 \label{eccb}
 \Gamma(D_p)  \ge i-\min({\bf A}).
\end{equation}
\end{proposition}

\begin{remark}
  For $p > 2$ the estimate (\ref{eccb}) is better than the estimate of the Proposition \ref{mhe}.
\end{remark}

\subsection{Estimates and calculation  of covering density }

At first give the estimation of the thinnest lattice covering of $D_p$ from above on the base of  results by Sas and Dovker (\cite{lek}).

\begin{proposition}
 The density  $\vartheta$ of any covering of  ${\mathbb R}^2$ by  Minkowski balls $D_p$ satisfies
  \begin{equation}
 \label{deu}
    \vartheta \le \frac{2\pi}{3\sqrt{3}}
\end{equation}

 \end{proposition}

 Let us now examine the curves of the family (\ref{sap}) from the point of view of the coverage densities they determine.
Here we do this only for one point on one curve.
 But even such example gives earlier unknown minimum density for Minkowski ball.

Take the curve from the family (\ref{sap} ) defined by the function
\begin{equation}
 \label{s2p}
  A(2,p) = 3(\tau + (2^p - 1)^\frac{1}{2 p})(1 + \tau^{p})^{-\frac{1}{p}}
  (1 + ((2^p - 1)^\frac{1}{2 p})^p)^{-\frac{1}{p}}. 
 \end{equation} 
 
 \begin{proposition}
 \label{emcd}
  The convex symmetric hexagon defined by the point $p =3$ on the curve (\ref{s2p}) gives a covering density $\vartheta(D_3) \approx 1.0567$ by the Minkowski ball $D_3$, the minimum for known densities of Minkowski balls at $1 < p < \infty$.
 \end{proposition}
 {\bf Proof}. \\
 ${\sigma}_{2,3} = (2^{3}-1)^\frac{1}{2\cdot3} \approx 1.3830;$\\
 $\tau_3 \approx 0.20406, 0 \le \tau \le \tau_3,$ for curve (\ref{s2p}) take $\tau = 0.1200;$\\
 $V(D_3) \approx 3.5200, \; \gamma_h(D_3) \approx 3.3310;  $ \\
 So $\vartheta(D_3) = V(D_3)/\gamma_h(D_3)   \approx 1.0567;$ \\
 
      We state here one conjecture and one problem which arise naturally from our work.
 
 \begin{conjecture*}     
 \label{bcc}
 The curve of maxima of covering constants increases from $p=1$ to $p=2$ and decreases from $p=2$ to $p=\infty$.     
\end{conjecture*}    
  
\begin{problem*}
\label{mbf}
Does the curve of maxima of covering constants belong to the family (\ref{sap})?
\end{problem*}

\section{Critical lattices, Minkowski domains and  their optimal packing}
   \label{secpmdccdomains}

   Let $D$ be a fixed bounded symmetric about origin  convex body ({\it centrally symmetric convex body} for short) with volume
$V(D)$.

\begin{proposition} \cite{Cassels,lek}.
\label{p1}
 If $D$ is symmetric about the origin and convex, then $2D$ is convex and symmetric
 about the origin.
\end{proposition}

For proves of results of this Section see \cite{gPacIIl}.

 \begin{corollary}
 \label{cor1}
 Let $m$ be integer $m \ge 0$ and $n$ be natural greater $m$. 
If $2^m D$ centrally symmetric convex body then $2^n D$ is again centrally symmetrc convex body.
\end{corollary}

\subsection{Domains.}

We consider the following classes of balls (see Section \ref{Minkowski balls}) and domains.

\begin{itemize}
  
\item  {\it  Minkowski domains}:    $2^m D_p$,  integer $m \ge 1$, for $1 \le p<2$;

\item {\it Davis domains}: $2^m D_p$,  integer $m \ge 1$, for $p_{0} > p \ge 2$;


\item {\it Chebyshev-Mordell domains}: $2^m D_p$,  integer $m \ge 1$,  for $ p \ge p_{0}$;

\end{itemize}

\begin{remark}
\label{r2}
Sometimes, when it comes to using $p$ values that include scopes of different domains, we will use the term Minkowski domains for definitions of different types of domains, specifying
  name when a specific $p$ value or a range of $p$ values containing a single type of domains is specified.
\end{remark}  

\begin{proposition}
\label{p7}
Let $m$ be  integer, $m \ge 1$.
If $\Lambda$ is the critical lattice of the 
ball $D_p$ than the sublattice $\Lambda_{2^m}$  of index $2^m$ is the critical lattice of the domain $2^{m-1} D_p$.
\end{proposition}

\begin{proposition} \cite{Cassels,lek}.
\label{p2m21m}
 The dencity of a $(2^{m-1} D_p, \Lambda)$-packing is equal to $V(2^{m-1} D_p)/d(\Lambda)$ and it is maximal if $\Lambda$ is critical for  $(2)^m D_p$.
\end{proposition}

\begin{theorem}
\label{opd}
The optimal lattice packing of  Minkowski and Chebyshev-Mordell domains $2^m D_p$ is realized with respect to the sublattices of the index $2^{m+1}$ of  critical lattices
  \[
  (-2^{-1/p},2^{-1/p}) \in \Lambda_{p}^{(1)}.
  \]   
The dencity of this optimal packing is given by the expression
\[
    V(2^m D_p)/\Delta^{(1)}_p (2^{m+1}D_p),
\]    
    where $V(2^m D_p)$ is the volume of the domainl $2^m D_p$ and $\Delta^{(1)}_p (2^{m+1} D_p)$ is the critical determinant of the domain $2^{m+1} D_p$.\\
The optimal lattice packing of  Davis domains $2^m D_p$ is realized with respect to the sublattices of the index $2^{m+1}$  of  critical lattices
  \[
  (1,0)\in\Lambda_{p}^{(0)},
  \]   
The dencity of this optimal packing is given by the expression 
\[
    V(2^m D_p )/\Delta^{(0)}_p(2^{m+1} D_p)
\]
 where $V(2^m D_p)$ is the volume of the domainl $2^m D_p$ and $\Delta^{(0)}_p (2^{m+1} D_p)$ is the critical determinant of the domain $2^{m+1} D_p$.\\
     \end{theorem}

\section{Inscribed and circumscribed hexagons of minimum areas}
\label{secichminarea}
For proves  see \cite{gPacIIl}.

Denote by $\Delta (2^m D_p)$ the critical determinant of the domain $2^m D_p$.

\begin{equation}
\label{cdd}
\Delta(2^m D_p) = \left\{
                   \begin{array}{lc}
     {\Delta^{(0)}_p}(2^m D_p), \; 1 \le p \le 2, \; p \ge p_{0},\\
    {\Delta^{(1)}_p}(2^m D_p), \;  2 \le p \le p_{0};\\
                     \end{array}
                       \right.
 \end{equation}
were $p_{0}$ is a real number that is defined unique by conditions
$\Delta(p_{0},\sigma_p) = \Delta(p_{0},1),  \;
2,57 < p_{0}  < 2,58, \; p_0  \approx 2.5725 $\\

Denoted by  $Ihma_{2^m D_p}$ the minimal area of  hexagons which inscibed in the domain $2^m D_p$  and have three pairs of points on the boundary of $2^m D_p$. 

\begin{theorem}
\label{ihma}
\[
  Ihma_{2^m D_p} = 3 \cdot \Delta (2^m D_p).
\]
 \end{theorem}
 
 Respectively denoted by  $Shma_{2^m D_p}$ the minimal area of  hexagons which circumscribed to the domain $2^m D_p$  and have three pairs of points on the boundary of $2^m D_p$. 

\begin{theorem}
\label{chma}
\[
        Shma_{2^m D_p}   = 4 \cdot \Delta (2^m D_p).
        \]
 \end{theorem} 

\section{On coverings by Minkowski domains in the plane}

\section{Direct systems,  direct limits and duality}
\label{secdsdld}
A more detailed description is given in \cite{gPacIIl}.

\begin{definition}
\label{dset}
 Let ${\mathbb N}_0$ be the set of natural numbers with zero.
A preoder $N$ on ${\mathbb N}_0$  is called a {\it directed set} if for each pair $k, m \in N$ there exists an $n \in N$ for which 
$k \le n$ and $m \le n$.
A subset  $N'$ is {\it cofinal} in  $N$ if, for each $m \in N$ their exists an $n \in N'$ such that $m \le n$.
\end{definition}

\begin{definition}
\label{dds}
A {\it direct} (or {\it inductive}) {\it system of sets} $\{X, \pi\}$ over a directed set $N$ is a function which attaches to each 
$m \in N$ a set $X^m$, and to each pair $m, n$  such that $m \le n$ in $N$, a map $\pi_m^n: X^m \to X^n$ such that, for each 
$m \in N$,  $\pi_m^m = Id$, and for $m \le n \le k$ in $N$, $\pi_n^k \pi_m^n = \pi_m^k$.
\end{definition}

\begin{remark}
We will consider direct systems of sets, topological spaces, groups and free ${\mathbb Z}$-modules.
\end{remark}

 \begin{proposition}
A directed set $N$ forms a category with elements are natural numbers $n, m, \ldots$  with zero and with morphisms
 $m \to n$ defined by the relation $m \le n$.~ A direct system over $N$ is a covariant functor from $N$ to the category of sets and maps, or to the category of topological spaces and continuous mappins, or to the category of groups  and homomorphisms or to the category of  ${\mathbb Z}$-modules and homomorphisms.
 \end{proposition}

\subsection{Direct systems of Minkowski domains and their limits}
\cite{gPacIIl}

\subsubsection{Direct systems}

Classes of direct systems of balls and Minkowski domains and classes of direct systems of their critical lattices are given in \cite{gPacIIl}.

\subsubsection{Direct limits}\cite{gPacIIl}

Let us give direct limits of these derect systems.
Let ${\mathbb Q}_2$ and ${\mathbb Z}_2$ be respectively the field of $2$-adic numbers and its ring of integers.

 \begin{proposition}
  $D^{dirlim}_p = \varinjlim  2^m D_p  \in ({\mathbb Q}_2 / {\mathbb Z}_2)  D_p =
   (\bigcup_m \frac{1}{2^m} {\mathbb Z}_2/{\mathbb Z}_2) D_p.$
 \end{proposition}
 
  {\bf  Proof}. Follow from properties of direct systems and their direct limits \cite{bur,stei}.

 \begin{proposition}
  $\Lambda_{p}^{dirlim} = \varinjlim  2^m \Lambda_{p}  \in ({\mathbb Q}_2 / {\mathbb Z}_2)  \Lambda_{p} = 
  (\bigcup_m \frac{1}{2^m} {\mathbb Z}_2/{\mathbb Z}_2) \Lambda_{p}.$
 \end{proposition}
 
 {\bf  Proof}. Follow from properties of direct systems and their direct limits \cite{bur,stei}.\\

 \section{About association schemes and spherical designs}

 Some results about invariant relations and group theory are given in papers \cite{gks,gkss}.
In this regard, let us recall the method of invariant relations and its origins.
 The connection between group theory and invariant relations of groups can be traced back to the works of E. Galois. General approaches to the theory of invariant relations were proposed by M. Krasner (see
  \cite{kra} and references therein). The method of invariant relations for symmetric groups was developed by Wielandt \cite{wie}.



 

 


 


\section{About Riemann-Roch theorem for algebraic curves,  lattices and arithmetic curves}
 For natural $p, n, p = 2n$ the curves 
 \[
 C_{2n}: x^{2n} + y^{2n} = 1
 \]
  are algebraic curves and after complexification their projective models have for $n > 1$  the genus  
\begin{equation}
 g = \frac{(2n - 1)(2n - 2)}{2}= (2n - 1)(n - 1).
\end{equation}

 \begin{theorem}
 (Riemann-Roch) 
\cite{har}. Let $D$ be a divisor on a curve $C_{2n}$ of genus $g$.
Then 
\begin{equation}
l(D) -  l(K - D) = \deg D + 1 - g
\end{equation}
 \end{theorem}
 
 \subsection{On absolute algebra and absolute algebraic geometry}
 The category-theoretic approach to number-theoretic and algebraic problems
 over finite fields ${\mathbb F}_q$ and their limit as $q$ tends to $1$ leads to 
interesting results and intriguing hypotheses (see \cite{tit,smi,ma,sou,dei,du,ccaa,ccag,ccsp,ccrrz}).

 In papers \cite{ccaa,ccag,ccsp,ccrrz,ccars,ccaar}, the authors construct absolute algebra and absolute algebraic geometry (affine case) and use it, in particular, to derive the Riemann-Roch theorem in the cases under study.  Let us recall very briefly the approaches of A. Connes and C. Consani to
 absolute algebra and absolute algebraic geometry.
  Let $k \ge 0$ be an integer and let $k_+$ be the set $k_+ = \{0, \ldots, k \}$, where $0$ is 
 the base point. Let  ${\Gamma}^{op}$ be the small, full subcategory of the category 
 ${\mathfrak{Fin}}_*$ of finite pointed sets, whose objects are  the sets $k_+$.
 A ${\Gamma}$-set is a functor ${\Gamma}^{op} \to {\mathfrak{Sets}}_*$ between these pointed categories.
 Here ${\mathfrak{Sets}}_*$ is the category of pointed sets.
   Absolute algebra connects with category   ${\Gamma}{\mathfrak{Sets}}_*$ 
 of ${\Gamma}$-sets. This is closed symmetric monoidal category \cite{dgm}.
 In the breakthrough paper \cite{ccag} authors present the affine case of absolute algebraic geometry as ``algebraic geometry for general Segal's $\Gamma$-rings''.
 Let  $ \mathfrak{Ring}$ be the category of commutative rings and $\mathfrak{Mo}$ be the category of commutative monoids. The category $\mathfrak{MR} = \mathfrak{Ring} \cup_{\beta, \beta^*}\mathfrak{Mo}$  is obtained by gluing togeter the categories $\mathfrak{Ring}$  and $\mathfrak{Mo}$ using the pair of adjoint functors $\beta$ and 
 ${\beta^*}$.
The main results in the paper are Theorem 3.1 on the category  $\mathfrak{MR}$ , Theorem 5.11 on the functoriality
 of ${\mathfrak {Spec}}(A)$,  Theorem 6.5 on the topos  ${\mathfrak {Spec}}(HR)$ in the case of a semiring $R$ and Theorem 7.10 on the functoriality of the structure sheaf $\mathcal O^{++}_A$.
 Let $\mathfrak S$ be a  $\mathbb S$-algebra by B. Dundas, T. Goodwillie, R. McCarthy \cite{dgm}.
Part (iii) of the Theorem 3.1 of the paper \cite{ccag} asserts that the category  $\mathfrak{MR}$ is a full subcategory of the category $\mathfrak S$ of  $\mathbb S$-algebras.
Now let $A$ be a commutative $\mathbb S$-algebra.
  A. Connes,  C. Consani   provide the basics of (affine) algebraic geometry over $\mathbb S$ by defining and investigating for ``$\mathbb S$-algebra $A$ its spectrum
${\mathfrak {Spec}}(A)$ understood as the pair of a Grothendieck site and a structure presheaf of $\mathbb S$-algebras,
 whose associated sheaf is a sheaf of $\mathbb S$-algebras''.
 The Grothendieck site is defined by authors as the pair $(C^{\infty}(A), J(A))$  where $C^{\infty}(A)$ is a small category and $J(A)$ is the Grothendieck topology on $C^{\infty}(A)$.
Let $R$ be a semiring and $HR$ be the corresponding $\mathbb S$-algebra. 
Theorem 6.5 compares the points of the topos  ${\mathfrak {Spec}}(HR)$ with the set of prime ideals of $R$,  and the site $(C^{\infty}(HR), J(HR))$ with Zariski spectrum $Spec(R)$.
The paper\cite{ccag}  is divided into nine sections and twenty one subsections. The first section is the Introduction and 
 is designed to acquire the reader with the basic ideas and results of the paper.
 Section 2 and 3 treat the connection between morphisms of $\mathbb S$-algebras and morphisms of $\mathfrak{MR}$ as well as functors that identify $\mathfrak{Mo}$, $\mathfrak{Ring}$  and $\mathfrak{MR}$ with full subcategories of the category $\mathfrak S$ of  $\mathbb S$-algebras.
 In section 4 a monoid $M$ of $\mathfrak{Mo}$ is studied, primarily for its importance to the case of the specific monoid $A(1_+)$ of an arbitrary $\mathbb S$-algebra $A$.
A small category $C(M)$ and a corresponding topos of presheaves  $\mathfrak{Root}(M)$ are presented and investigated.
Sections 5 - 7  are the heart of the paper. 
Section 5  gives  ``the basis to develop algebraic geometry over  $\mathbb S$. The main construction, in the affine case, is to associate to a (commutative)  $\mathbb S$-algebra $A$ its spectrum
$\mathfrak{Spec}(A)$  understood as the pair of a Grothendieck site and a structure presheaf of $\mathbb S$-algebras, whose associated sheaf is a sheaf of $\mathbb S$-algebras''. Section 6 presents examples of $\mathfrak{Spec}(A)$.
Section 7 treats the structurel sheaf on $\mathfrak{Spec}(A)$. Results of the section are applicable to the adele class space of the first author \cite{ctf}.
In section 8 authors ``discuss the relations between the above developments and the theory of Toën B., 
 Vaqui\`e'' \cite{tva}.

  \subsection{On Riemann-Roch}
  Let $D$ be the Arakelov divisor (replete divisor by J. Neukirch \cite{neu})  on $\overline{Spec \; {\mathbb Z}}$.
   The main result of the paper \cite{ccrrz} is 
   \begin{theorem}
  \cite{ccrrz}. (Riemann-Roch formula for $D$)
  \begin{equation}
    dim_{\mathbb S} H^0(D) - dim_{\mathbb S} H^1(D) =  \lceil \deg_2 D  \rceil^\prime,
    \end{equation}
   \end{theorem}
 where $\deg_2 := \deg / \log 2$, $\lceil x \rceil^\prime$  ``denotes the right continuous function which agrees with the function $ceiling(x)$ for $x > 0$ non-integer, and with $-ceiling(-x)$ for $x < 0$ non-integer'', and ${\mathbb S}$ is the absolute base.
This formula, obtained by the authors while working on the absolute base itself, improves the most spectacular result (the Riemann-Roch formula for $\overline{Spec\; {\mathbb Z}}$) obtained while working on the spherical extension ${\mathbb S}[\pm 1]$ of the paper \cite{ccsp}.

\section{On algebraic dynamics related to Minkowski's hypothesis about the critical determinant of
$|x|^p + |y|^p  < 1 \; p \ge 1$ domain}

 \subsection{On algebraic dynamics and resurgence on Minkowski moduli spaces \cite{gefr,adr,glkh}} 
For ease of reading, let us recall the notation.
Let 
$$    
 |x|^p + |y|^p \le 1, \; p > 1,                     
$$
be the solid Minkowski tube
with the surface boundary ${\mathcal MT}$.
\begin{remark}
 For natural $p, d, p = 2d$ the Minkowski curves 
\begin{equation}
\label{mac} 
  {\mathbb S}^{1}_{2d} = C_{2d}: x^{2d} + y^{2d} = 1
\end{equation}
  are algebraic curves. They define algebraic layers of the moduli space ${\mathcal MT}$ and the bundle of their Jacobians 
after complexification of the corresponding bundles.
\end{remark}

The Minkowski-Cohn moduli  space ${\mathcal M}$ of admissible lattices of Minkowski balls has the form
 
 \begin{equation}
 \label{mcms}
 \Delta(p,\sigma) = (\tau + \sigma)(1 + \tau^{p})^{-\frac{1}{p}}
  (1 + \sigma^p)^{-\frac{1}{p}}, 
 \end{equation}
in the domain
\begin{equation}
\label{msd}
  {\mathcal M}: \; \infty > p > 1, \; 1 \leq \sigma \leq \sigma_{p} =
 (2^p - 1)^{\frac{1}{p}},
  \end{equation}
of the $ \{p,\sigma\} $-plane, where $\sigma$ is some real parameter.
\subsubsection{$l$-adic representations and Tate modules}
Let $K$ be a field and $\overline K$ its separate closure,
$$ E_n = \{P \in E(\overline K)|nP = 0\}$$
 the group of points
of elliptic curve  
$E({\overline{K}})$ 
order dividing $n$. 

When 
$char \; K$ does not divide $n$ then $E_n$ is a free 
${\mathbb{Z}}/n{\mathbb{Z}}$-module of rank $2.$
Let $l$ be prime, $l \neq char \; K$. The projective limit $T_l(E)$ of the projective system of modules $E_{l^m}$
is free ${\mathbb{Z}}_l$-adic Tate module of rank 2.

Let $V_l(E) = T_l(E) \otimes_{{\mathbb{Z}}_l}
  {\mathbb{Q}}_l 
$.

 Galois group $Gal({\overline{K}}/K)$ 
acts on all $E_{l^m}$, so there is the natural
continuous representation ($l$-adic representation)

$$ 
\rho_{E,l}: Gal({\overline{K}}/K) \to Aut(T_l(E))  \subseteq Aut(V_l(E)).
$$

For Jacobians of Minkowski curves (\ref{mac}) we have Tate modules of dimension $2g$ and corresponding $l$-adic representations.

\begin{proposition}
These representations define $l$-adic Lie groups that act on on corresponding $l$-adic spaces and generate corresponding dynamics.
\end{proposition} 

\subsubsection{Infinitesimal deformations of Lie groups}

The corresponding infinitesimal deformations are associated with formal schemes and formal groups \cite{glp,glkh}.\\

\subsubsection{Differential-algebraic equations \cite{kkg}} 
Roughly speaking a differential-algebraic equation (DAE)  is a differential equation on algebraic variety.

The system of DAE has the form
\begin{equation}
 F=F(dx/dt, x, t) = 0, F = (F_1, \ldots, F_n)^T,  
\end{equation}
$
                                              x(t) = (x_1(t), \ldots, x_n(t)) \in{\mathbb R}^n, t \in {\mathbb R},
$

One of the algebraic-geometric models of DAE is the Gauss-Manin connection on an algebraic variety\\

\subsubsection{About the  resurgence \cite{etr,pon}}

This term have introduced (defined) by J. Ecalle in 1998 year \cite{etr}.

Consider Euler's equation
\begin{equation}
\label{ee}
 t^2( dy/dt) = t -y
\end{equation}

By a
book  on ODE
\cite{pon}
 we try to solve (\ref{ee}) by

$
y = \sum_{m=0}^{\infty} a_m t^m
$

We obtaine

$
 y = \sum_{m=0}^{\infty}(-1)^m m! t^{m+1}
$

Unfortunately this series divergences  for all $t \neq  0$.

But if we consider homogeneous equation 
$  t^2( dy/dt) =  -y  $
and solve it, we have 

$ y = ae^{1/t}.$

The solution
$$
   y = \sum_{m=0}^{\infty}(-1)^m m! t^{m+1} + ae^{1/t}.
$$ 
is the resurgence solution of  (5).

I will call the resurgence the Tauberian resurgence.

There exists also the Abelian resurgence. 

We define on (\ref{mcms}) by resurgence semi complex vector bundle, elliptic curves bundle, bundle of  Epstein zeta functions, and investigate their algebraic and dynamical properties \cite{gefms}. .

Epstein zeta function

Let $Q(s)$ be a positive integral quadratic form.

$$ \zeta_Q(s) = \sum_{(m,n) \neq (0,0)} \frac{1}{Q(m,n)^s}  $$

Let $s = \sigma + it.$

Laplace transform of a real function $f(u).$

$$
   \int_0^{\infty} e^{-su}f(u) du
$$

 \subsection{Algebraic non-Archimedean Minkowski curves and their analytification}

As we have mentioned earlier, for natural $p, n, p = 2n$ the curves
 \begin{equation}
\label{mincur}
C_{2n}: x^{2n} + y^{2n} = 1
 \end{equation}
 are algebraic curves. These curves are defined over real numbers.

For arithmetic and algebraic applications, equations are  also defined over non-Archimedean fields and rings.

 Let a curve (\ref{mincur}) be defined over the completion with respect to the non-trivial absolute value $|\cdot|_K$ of an algebraically closed non-Archimedean local field $K$.\\
For natural $n$, starting with $n = 1$, one can define  connected smooth strictly $K$-analytic curves $C_{2n}^{an} = {\mathcal C}_{2n}$ in the sense by Berkovich \cite{ber}

  \section{On the study of activity measures of dynamic systems over non-Archimedean fields}
Consider an analytic family of rational functions of degree $d$ parametrized by a connected smooth strictly $K$-analytic curve $V$ in the sense by Bercovich \cite{ber}.
Here $\mathcal K$ is an algebraically closed field with a complete, nontrivial and non-archimedean valuation. Bifurcation phenomena of a dynamic system are associated with the properties of points in the state space of the system.
The approach by DeMarco \cite{dem} in complex dynamics gives a new condition for passivity of the pair $(f, c)$. Applications of Berkovich spaces to complex dynamics are given in the papers by M, Baker, L. Demarco \cite{bdmp} and by X. Yuan, S.-Zhang [Calabi theorem and algebraic dynamics, url={https://api.semanticscholar.org/CorpusID:15166330 (2009)].
Recall that the affine line ${\mathbb A}^{1,an}_K$ (or ${\mathbb A}^{1,an}$) of Berkovich is the set of multiplicative seminorms on $K[T]$ that induce the given absolute value $|·|$ on $K$ and the Berkovich projective line is ${\mathbb P}^{1,an} = {\mathbb A}^{1,an} \bigcup (\infty)$. On seminorms and non-classical points in the non-Archimedean case, see the monograph by Berkovich, and also, for example, the article by Baldassarri F. [Non-Archimedean gauge seminorms arXiv:1611.04753v1 [math.AG] (2016; Zbl)], and for Archimedean case, see, e.g., the book by Narici L., Bekenstein E. [Topological Vector Spaces, Boca Raton, FL: CRC Press (2011; Zbl)].
The author of the paper \cite{iro} investigates stability (or passivity) notions for the family of rational functions over $V$ and constructs the activity measures (of a critical point) of a family of rational functions.  In the study of stability, the author \cite{iro} uses the non-Archimedean version of Montel's theorem proposed by \cite{fkt} and in the study of activity measures, the potential theory on Berkovich spaces proposed by M. Baker and R. Roumely \cite{brpt}; A. Thuillier [PhD diss., University of Rennes 1 (2005; Zbl)]. For each $t \in V$ let ${\mathcal H}(t)$ be the complete residue field at $t$.
The author's \cite{iro} analytic family is an analytic morphism $f : V \times_K {\mathbb P}^{1,an} \to V \times_K {\mathbb P}^{1,an}$ preserving fibers, with an additional assumption that for each $t \in V$ , the morphism $f_t := f(t, .) : {\mathcal H}(t) \times_K {\mathbb P}^{1,an} \to {\mathcal H}(t) \times_K {\mathbb P}^{1,an}$ is a rational function over ${\mathcal H}(t) \times_K {\mathbb P}^{1,an} \simeq {\mathbb P}^{1,an}_{{\mathcal H}(t)}$ of degree $d$. The author fixes an analytic morphism $c: V \to {\mathbb P}^{1,an}$ and calls it a marked point. `The purpose of this paper \cite{iro} is to develop a stability theory for such a pair $(f, c)$ i.e. to study the asymptotic properties of $\{c_n(t)\}^{\infty}_{n=0}$`. Here $c_n(t)(=f^n_t(c(t)))$ is the composition $\pi_2 \circ f^n \circ (id \times c) : V \to {\mathbb P}^{1,an}$, where $\pi_2 : V \times_K {\mathbb P}^{1,an} \to {\mathbb P}^{1,an}$ is the second projection. Theorem 1.1. \cite{iro} Let $K, V, f, c$ and $c_n$ be as above. Then, there exists a unique positive Radon measure $\mu_{f,c}$ on $V$ such that for any non-classical point $\zeta \in V$ , we have the weak convergence of measures $\lim_{n \to \infty}\frac{1}{d^n}c^{*}_n \delta_{\zeta} = \mu_{f,c}$, when $\delta_{\zeta}$ is the Dirac mass at ${\zeta}$. Author calls $\mu_{f,c}$ the activity measure of the pair $(f,c)$. The proof given in this paper \cite{iro} uses potential theory (with corresponding Laplace operator) over Berkovich curve and `a certain equidistribution result, which also implies the uniqueness in Theorem 1.1 and the independence of the choice of the coordinate on ${\mathbb P}^{1,an}$'. The second main theorem concerns the support of the activity measure of a family of polynomials. Let the residue characteristic of $K$ is $0$ or greater than $d$. For a certain algebraic family of monic polynomials over $V = {\mathbb A}^{1,an}_K$ let ${\mathcal M}_{f,c} = \{t \in V| \{c_n(t)\} \mbox{is bounded}\}$ be the boundedness locus for the pair $(f,c)$ with $c$ critical. Let ${\mathcal M}_{f,c}$ be non-empty and bounded. Theorem 4.2 \cite{iro}. Let $(f, c)$ be a pair consisting of a family $f$ of monic polynomials and a marked critical point $c$. Then, the activity measure $\mu_{f,c}$ is proportional to the equilibrium measure of the set ${\mathcal M}_{f,c}$ with respect to $\infty$. In particular, the support of the activity measure coincides with the boundary of ${\mathcal M}_{f,c}$. In the case of the sutisfaction of conditions of Theorem 1.2 \cite{iro} under assumption, that the Mandelbrot set ${\mathcal M}_{f,c}$ is non-empty and bounded, two notions of activity are equivalent by author`s Proposition 5.2 \cite{iro}.


\section{On the theory of matroids}
\subsection{On definition of matroids}
There are several equivalent definitions of the concept of matroids.
Below we give definitions that can be used not only for finite matroids.
 We recall three of them: by independent sets, by circuits, and by flats. Below we will denote the cardinality of a (finite) set $S$ by 
$\#S$.

 (IM) Definition by independent sets: A matroid $M = (E,{\mathcal I})$ is a pair consists of a finite set $E$ and a collection ${\mathcal I}$ of subsets of $E$, satisfying the following axioms:\\
(IM-1) ${\emptyset} \in {\mathcal I}$, ${\mathcal I} \ne {\emptyset} $. \\
(IM-2)  if $I \in {\mathcal I}$ and $J  \subseteq I$ then $J \in {\mathcal I}$.\\
(IM-3)  if $I, J \in {\mathcal I}$ and  $\#I > \#J$, then there exists $i \in I - J$ such that $J \cup i \in {\mathcal I}$.    

Any maximal independent set of $M$ is called the basis of $M$.
All maximal independent sets (all bases) of $M$ have the same 
cardinality. The cardinality (the size) of a basis is called the rank $r = r(M)$ of $M$.

\subsubsection{Definition by the rank function} 
The rank function can be defined for any subset $S \subseteq E$: $r(S)$ is the common cardinality of the sets of $\mathcal I$ which are maximal in $S$. So $r$ is a rank function on the subsets of $ E$ having the properties: \\
(R-1) $ 0  \leq r(S)  \leq \#S $;\\
(R-2) If $S \subseteq T$ then $r(S) \leq r(T)$;\\
(R-3) $r(S \cup T) + r(S \cap T) \leq r(S) + r(T)$;\\

The rank function with properties (R-1) - (R-3) (and with known values on all subsets of $E$) defines the matroid on the set $E$: the subsets $S$ of $E$ for which $r(S) = \#S$ 
are the independent  sets of a matroid.

Let $\mathcal B$ be the set of all bases of $M$.\\
(CM) Definition by circuits: A matroid $M = (E,{\mathcal C})$ is a pair consists of a finite set $E$ and 
 a collection $\mathcal C$ of non-empty incomparable subsets of $E$, called circuits, satisfying
  the axiom:\\
(CM) if $C_1$ and $C_2$ are distinct circuits
and $c \in C_1 \bigcap C_2$, then   $(C_1 \bigcup C_2) - \{c\}$ contains a circuit. \\

(FM) Definition by flats:   A matroid $M = (E,{\mathcal F})$ is a pair consists of a finite set $E$ and 
 a collection $\mathcal F$  of subsets   of $E$, called flats, satisfying the following axioms:\\
 (FM-1) If $F_1$ and $F_2$ are flats of $M$, then their intersection is a flats of $M$. \\
  (FM-2) If $F \in {\mathcal F}$,  then any element of $E \setminus F$ is contained in exactly one flat of $M$ covering $F$.\\
   (FM-3)   The ground set $E$ is the flat of $M$.

   \begin{remark}
     $\varnothing \in {\mathcal F}$. \\
    It is say that a flat $G$ of $M$ is a cover of $F$ if $G$ is minimal among the flats of $M$ properly containing $F$.
   \end{remark}
   
    \begin{definition}
  Two matroids $M = (E_1, {\mathcal I_1})$ and  $N = (E_2,{\mathcal I_2})$ are isomorphic if  there is a relabeling bijection $\alpha: E_1 \to E_2$ that maps $I_1$ to $I_2$.
  \end{definition}
   
   \subsection{Matroids and some topological notions} 
    In definition (IM) ${\mathcal I}$ is a nonempty (abstract) simplicial complex. Its geometric realization is a polytope.
      
{\it {Closures and  matroids}}  
 \cite{wiki}.\\
Let 
 ${\displaystyle M} = ({\displaystyle E}, {\mathcal I})$ be a matroid on a finite set
${\displaystyle E}$, with rank function 
${\displaystyle r} $.

Let $\overline S$ be the set of $x \in E$ such that $r(S \cup \{x\}) = r(S)$. This defines the closure operetor 
$cl: S \to \overline S$.  More formally:

the closure or span ${\displaystyle \operatorname {cl} (S)}$ of a subset 
${\displaystyle S} $ of
${\displaystyle E}$ is the set

${\displaystyle \operatorname {cl} (S)={\Bigl \{}\ x\in E\mid r(S)=r{\bigl (}S\cup \{x\}{\bigr )}{\Bigr \}}}$,

so we have: 

${\displaystyle \operatorname {cl} :{\mathcal {P}}(E)\mapsto {\mathcal {P}}(E)}$ 

where 

${\displaystyle {\mathcal {P}}}$ denotes the power set, with the following properties:

(C1) For all subsets 
${\displaystyle X}$ of 
${\displaystyle E}$, 

${\displaystyle X\subseteq \operatorname {cl} (X)}$.

(C2) For all subsets 
${\displaystyle X}$ of 
${\displaystyle E}$,

${\displaystyle \operatorname {cl} (X)=\operatorname {cl} \left(\operatorname {cl} \left(X\right)\right)}$.

(C3) For all subsets 
${\displaystyle X}$ and 
${\displaystyle Y}$ of 
${\displaystyle E}$ with 

${\displaystyle X  \subseteq Y}$

${\displaystyle \operatorname {cl} (X)\subseteq \operatorname {cl} (Y)}$.

(C4) For all elements 
${\displaystyle a}$ and 
${\displaystyle b}$ from 
${\displaystyle E}$

 and all subsets 
${\displaystyle Y}$ of 
${\displaystyle E}$, if 

${\displaystyle a\in \operatorname {cl} (Y\cup \{b\})\smallsetminus \operatorname {cl} (Y)}$

 then 

${\displaystyle b\in \operatorname {cl} (Y\cup \{a\})\smallsetminus \operatorname {cl} (Y)}$.

The first three of these properties are the defining properties of a closure operator.

${\displaystyle \operatorname {cl} :{\mathcal {P}}(E)\to {\mathcal {P}}(E)}$
 that obeys these properties determines a matroid.

{\it {Closures and Flats}}\\ 
A set whose closure equals itself is said to be closed.
\begin{remark}
 A closure set is a flat (or subspace of the matroid) [4].
A set is closed if it is maximal for its rank, meaning that the addition of any other element to the set would increase the rank.
\end{remark}

\begin{remark}
  By the defonition of the closure operation we have (the closed sets of a matroid are characterized by a covering partition property):
 \\
(F1) The whole point set 
${\displaystyle E}$ is closed.

(F2) If 
${\displaystyle S}$ and 
${\displaystyle T}$ are flats, then 
${\displaystyle S\cap T}$ is a flat.

(F3) If 
${\displaystyle S}$ is a flat, then each element of 
${\displaystyle E\smallsetminus S}$ 

is in precisely one of the flats 
${\displaystyle T}$ that cover 
${\displaystyle S}$ 

(meaning that 
${\displaystyle T}$ properly contains 
${\displaystyle S}$, but there is no flat 
${\displaystyle U}$ between 
${\displaystyle S}$ and 
${\displaystyle T}$)
\end{remark}

 \subsection{Matroid classes}
    \begin{definition}
  Two matroids $M = (E_1, {\mathcal I_1})$ and  $N = (E_2,{\mathcal I_2})$ are isomorphic if  there is a relabeling bijection $\alpha: E_1 \to E_2$ that maps ${\mathcal I_1}$ to ${\mathcal I_2}$.
  \end{definition}
   
      \subsubsection{Uniform matroids}
  Let    $M = (E,{\mathcal I}), \: \#E = n,$ be a matroid and let all bases of it consist of $k \le n$ element subsets, {\it id est}  $\forall B \in {\mathcal B}, \#B = k$. Such matroid is called uniform and is denoted by $U_n^k$ (or by $U_{k,n}$).\\
  
  \subsubsection{Free matroids}
   Let    $M = (E,{\mathcal I}), \: \#E = n,$ be a matroid and let all bases of it consist of $ n$ element subsets.
 Respectively all independent sets are all subsets of $E$. It is a special case  $U^n_n$ of a uniform matroid.
 
   \subsubsection{Simple matroids}
   A matroid is called simple if it has no circuits consisting of 1 or 2 elements.
  
\subsubsection{Linear matroids} Let $K$ be a field. Let    $M = (E,{\mathcal I})$ be a matroid such that the set $E$ is the set of vectors of a vector space $V_K$ over $K$, and let ${\mathcal I}$ be the collection of liniarly independent subsets of the set $E$.
The matroid $M = (E,{\mathcal I})$ is called the linear matroid.

\subsubsection{Graphical matroids} Let $E$ be the set of edges of a graph $G$ and let ${\mathcal I}$ be the collection of forests of 
$G$. Then the  matroid $M = (E,{\mathcal I})$ is a graphical matroid.
The class 
${\displaystyle {\mathcal {L}}(M)}$ 
 of all flats, partially ordered by set inclusion, forms a matroid lattice.
 
  Conversely, every matroid lattice 
${\displaystyle L}$ forms a matroid over its set 
${\displaystyle E}$ 
 of atoms under the following closure operator: for a set 
${\displaystyle S}$ of atoms with join 
${\displaystyle \bigvee S}$.

   \subsubsection{Uniform matroids}
  Let    $M = (E,{\mathcal I}), \: \#E = n,$ be a matroid and let all bases of it consist of $k \le n$ element subsets, id est  $\forall B \in {\mathcal I}, \#B = k$. Such matroid is called uniform and is denoted by $U_{k,n}$.\\
  
  \subsubsection{Linear matroids}
  Let $K$ be a field. Let $V$ be a vector space over $K$ and $E$ be a set of vectors of $V$. Let
  ${\mathcal I}$ be the collection of linear independent subsets of $E$. Then $M = (E,{\mathcal I})$ is a 
   linear matroid over $K$.
   
   Let $S$  be a partially ordered set and $\omega$
 be a symbol with $\omega  \notin S$.   Let $K$ be  a field and $V_s$  a $K$-vector space.
  An $S$ - space is of the form
 $V=(V_{\omega}, V_s)_{s\in S}$, where the $V_s$  are subspaces of
 the $K$- space $V_{\omega}$  for $s \in S$, such that $s\le s'$ 
 implies $V_s \subset V_{s'}$.

\subsection{ Operations with matroids}

 \subsubsection{Restriction}
 Let $S \subseteq E$ be a subset of the matroid $M = (E,{\mathcal I})$. The restriction $M\mid S$ 
 is the matroid on the ground set $S$ with independent sets ${\mathcal I}\mid S = \{{\mathcal I} \subseteq S, I \in {\mathcal I}\}$.

 \subsubsection{Contraction}
  Let $S \subseteq E$ be a subset of the matroid $M = (E,{\mathcal I})$. The contraction $M / S$ 
 is the matroid on the ground set $E - S$  with independent sets 
 ${\mathcal I} / S = \{{\mathcal I} \subseteq E \setminus S, I \cup I_S \in {\mathcal I}\}$ 
 for any maximal independent subset $I_S$ of $S$.
 
    \begin{remark} 
    Sometimes these two operations (restriction and   contraction) are called reduction (of a matroid).
        \end{remark}
    
   \begin{definition}  
It is said that a matroid $N$ is a minor of a matroid $M$ if it can be constructed
from $M$ by restriction and contraction operations (by a  reduction  the matriod). 
\end{definition}
 
 

\subsection{Characteristic polynomial of a matroid}

Let $M$ be a linear matroid with the rank function $r$. The characteristic polynomial of $M$ is
\begin{equation}
\label{chp}
 \chi_M(q) = \sum_{S \subseteq E} (-1)^{\#S} q^{r(M) - r(S)}.
\end{equation}
Let $\chi_M(q) = w_0 q^r - w_1 q^{r-1} + \cdots + (-1)^r w_r. $ The numbers $w_i, 0 \leq i \leq r$ are called Whitney numbers.
 \;

\section{On elements of the history of matroids}
   Matroids (as finite matroids) were defined independently by Whitney \cite{whi} and Nakasawa \cite{nak}.
G. Birkhoff noted the connection between matroids and 
geometric lattices.
S. MacLane noted that the properties of transcendence bases
of field extensions satisfy the axioms of a matroids.  
W. Tutte \cite{tut} explored the fundamental connections between matroids and graphs.
R. Rado \cite{rado} discovered the connection between matroids and combinatorics based on the matroid analogy of P. Hall's theorem on representatives of sets.
H. Crapo and G. Rota \cite{crarot} developed foundations of combinatorial theory on the base of geometric lattices.\\
The history of matroid theory up to 1986 year is described in the anthology by J. P. S. Kung \cite{kung}.
Below we very briefly present elements of the history of matroids based on reviews of some papers by J. Oxley \cite{oxl}, T. Zaslavsky \cite{zas,zas1} and   D. Mayhew, B. Oporowski, J. Oxley, G. Whittle \cite{moow}.

\subsection{On interplay between graphs and matroids  by J. Oxley}
 The main goal of paper \cite{oxl} is to give a pedagogical introduction to
 the subject of interplay between graphs and matroids. It is primary
 intended for graduate students and specialists who are interested 
in graphs, linear algebra, matroids, coding theory and their applications.
 This paper should be accessible for any reader with some background in 
graph theory and linear algebra. The basic ideas and techniques are 
treated in detail, including examples. This is a renew version of some
 chapters of the author’s book \cite{oxl1}.
 Also it includes  results (until 2002 year) on matroids.
The organization of the
 paper is as follows 1. Introduction. 2. The bare facts about matroids.
 3. Connectedness, 2-connectedness, and unavoidable structures. 
4. Some proofs. 5. Removable circuits. 6. Removing circuits from 
3-connected matroids. 7. Minors and infinite antichains. 
8. Branch-width and infinite antichains. 9. The implication of large
 branch-width. 10. Some proof outlines. 11. Unavoidability revisited.
The material of the paper is outlined in the titles of sections and
 informally discussed in the introduction. Section 2 gives an admirable
 straightforward introduction to the basics of matroid theory. 
“Section 3 begins by showing how 2-connectedness for graphs extends 
naturally to matroids. It then indicates how the number of edges in
 a 2-connected loopless graph can be bounded in terms of the 
circumference and the size of a largest bond. The main result of the
 section extends this graph result to matroids.” Section 4 outlines
 the proofs of the main results from Section 3. In Sections 5-6 the 
author describes verious extensions and non-extensions for 2-connected
 graphs and matroids in a 1974 theorem of W. Mader [Abh. Math. Semin.
 Univ. Hamb. 42, 187-204 (1974; Zbl 0266.05111)]. The graph minors 
project by N. Robertson and P. D. Seymour [Surveys in combinatorics 1985,
 Lond. Math. Soc. Lect. Note Ser. 103, 153-171 (1985; Zbl 0568.05025)] 
and a longstanding matroid conjecture of G.-C. Rota [Actes Congr. 
Internat. Math. 1970, 229-233 (1971; Zbl 0362.05044)] motivate in 
particular Sections 7-10. The paper ends with a comprehensive list of
 references, arranged alphabetically.

  \subsection{On Supersolvable Frame-matroid and Graphic-lift Lattices}
Zaslavsky in the papers \cite{zas,zas1} investigates sufficient conditions (supersolvability) for a complete integral factorization
 of the characteristic polynomials of geometric lattices that correspond to arrangements of hyperplanes 
and finite matroids. Let $L$ be  a geometric lattice. A flat $x$ of a lattice  is modular \cite{smi}, if for all $y, z \in L$
with $z \le y$ we have  $z \vee(x\wedge y) = (z\vee x)\wedge y$. Recall that a lattice $L$  of rank $n$
 is supersolvable \cite{smi} if there exists a saturated chain $0 = x_0 < x_1 < \cdots <x_{n-1}<x_n=1$
 with rank $(x_i)=i$  such that every flat $x_i$  is modular.    
  The paper \cite{zas1} includes interesting results on the supersolvability of lattices. In part
 it is a continuation of the article by R. Stanley \cite{sta1}.
The paper \cite{zas1} also includes a concisely-written survey of recent results on biased and gain 
graphs, biased and lift matroids, coloring and polynomials (Section 1). The main results of the paper 
are theorems about supersolvability of two types of geometric lattices: (i) lattices whose matroid is 
a frame; (ii) lattices defined by a graphic lift matroid. The author also gives the geometric 
interpretation of his theorems in terms of arrangements in $F^n$
, where $F$ 
 is a skew field. Arrangements of hyperplanes connect the integral representation of hypergeometric 
functions and Orlik-Solomon algebras with conformal field theory, representation theory, algebraic K
-theory and algebraic geometry. A clear and informative overview of these connections can be found 
in lectures notes by A. Varchenko \cite{var} and P. Orlik and
 H. Terao \cite{orte}.
The organization of the paper \cite{zas1} is as follows: Introduction. 1. Biased graphs, matroids,
 etc. 2. Frames. 3. Graphic lifts. 4. Examples. 5. Comments and questions.
The author calls the matroid of a geometric lattice (of finite rank) such that it has a basis, or 
can be extended to have one, and such that each point lies on a line generated by a pair of basis 
elements a (finite rank) frame matroid. The author defines a biased graph $\Omega = (V, E, B)$ as a graph 
$||\Omega|| = (V, E)$, not necessarily finite, together with a linear subclass $B = B(\Omega )$ of its polygons.
 The bias matroid $G(\Omega)$  is the matroid on the set of edges of a biased graph $\Omega$ whose circuits are 
the bias circuits.

Theorem 2.2. Let  $\Omega$
 be a simply biased graph of finite order. $G(\Omega)$
 is supersolvable if and only if each connected component of $\Omega$
 either: \\
(i) has a bias-simplicial vertex ordering; or \\
(ii) is a simplicial extension of one of \\
(a) $(mK_2,0)$,
 when $m \ge 2$,
 or \\
(b)  $⟨\pm K3⟩$, 
 or \\
 (c) $⟨\Sigma⟩$
 for $\Sigma = +\Gamma \cup - S_k$,
 where $\Gamma$
 is a chordal simple graph of finite order, $S_k$
 is a $k-$edge star whose vertex set lies in $V(\Gamma)$,
 and the noncentral vertices of $S_k$
 are a clique in $\Gamma$
. Furthermore, the bias-simplicial vertex ordering (or simplicial extension) can be chosen so that any 
desired bias-simplicial (or simplicial) vertex is the last vertex.

Section 3 contains the main theorems. For a biased graph $\Omega$
 the author defines the extended lift matroid $L_0(\Omega)$
 and the lift matroid $L(\Omega)$.\\
Theorem 3.2. Let $\Omega$ 
be a simply biased, connected link graph of finite order.\\
 (A) $L_0(\Omega)$ 
 is supersolvable if and only if $\Omega$
 has a link-simplicial vertex ordering.\\
(B) $L(\Omega)$ 
 is supersolvable if and only if:\\
 (i) $\Omega$
 is balanced and $||\Omega||$
 is chordal; or\\
 (ii) $\Omega$
 is as in Theorem 2.2 (ii). \\
Furthermore, the link-simplicial vertex ordering (or simplicial
 extension) can be chosen so that any desired bias-simplicial (or simplicial) vertex is last.\\

The author applies his results to three kinds of examples and to an extension of a theorem of
 P. H. Edelman and V. Reiner \cite{edre}.

\subsection{On excluded minors for the class of matroids that are binary or ternary}

The purpose of the paper \cite{moow}
 is to give the classification of the excluded minors for classes of matroids. \\

Theorem 1.1.: The excluded minors for the class of matroids that are binary or ternary are 
$U_{2,5}, U_{3,5}, U_{2,4}\oplus F_{7}, U_{2,4}\oplus F_{7}^{*}, U_{2,4}{\oplus}_{2} F_{7}, U_{2,4}{\oplus}_{2} F_{7}^{*}$, 
and the unique matroids obtained by relaxing a circuit-hyperplane in either $AG(3,2)$ or $T_{12}$.

The technique is matroid-theoretic and graph-theoretic, involving results from the book and paper 
by J. Oxley \cite{oxl1,oxl}  and results on the structure of almost-regular 
matroids by K. Truemper \cite{tru}, and leads to a finite task.
A number propositions relative to excluded minors that have low rank, corank or connectivity as well as 
excluded minors with at most nine elements are established.

These results show that authors of the paper \cite{moow} can restrict their attention to 3-connected excluded 
minors with rank and corank at least four and with at least ten elements and are then utilized to prove 
Theorem 5.1 on producing of such excluded minors by relaxing a circuit-hyperplane in a binary matroid.

Let $\mathcal M$ be a class of matroids that are either binary or ternary. The authors prove Theorem 6.1 and 
Proposition 6.2 concerning a connection between the excluded minors of $\mathcal M$ and Truemper`s class of 
almost-regular matroids. After the reduction to a finite list of excluded minors authors conclude with a proof of 
their principal result.

A number of examples and figures are provided to illustrate the basic theorems.
There is a bibliography of thirty one  items.
 
\section{Matroids and Minkowski  spheres} 
 Recall the definition of Minkowski's spheres 
${\mathbb S}^{n-1}_c$ (see equation (\ref{mt}))$$
{{\mathbb S}^{n-1}_c}: \;   |x_1|^c + \cdots + |x_n|^c  = 1, \;  c > 1.
$$
Here we consider the case of natural $n, \: 2 \le n \le 4$.
Consider matroids that may correspond these Minkowski's spheres.\\

One of the possible classes of such matroids will correspond to 
Euclidean spheres of the above dimensions, having a currently known
 quantity of kissing numbers.
 
 Let us recall in this connection the definition of one class of one-dimensional spheres.
 See Definition \ref{dc}.
 
The class of {\it Chebyshev-Mordell} spheres ${\mathbb S}^{1}_c$ is defined  for $ c > p_{0}$;
here $p_{0}$ is a real number that is  uniquely determined by the conditions
$\Delta(p_{0},\sigma_p) = \Delta(p_{0},1),  \;
2,57 < p_{0}  < 2,58, \; p_0  \approx 2.5725 $.

\begin{proposition} \cite{CS, mus}
Let $n = 2$. The  kissing number of  the Euclidean sphere ${\mathbb S}^{1}$ is equal 6 (six).\\
Let $n = 3$.The  kissing number of  the Euclidean sphere ${\mathbb S}^{2}$ is equal 12 (twelve).\\
Let $n = 4$. The  kissing number of  the Euclidean sphere ${\mathbb S}^{3}$ is equal 24 (twenty four).
 \end{proposition}
 
 Similar to the problem on the maximum possible kissing number 
for $(n - 1)$-dimensional  Euclidean spheres in $n $-dimensional
 Euclidean space, one can pose the following problem:
\begin{problem}
 What is the maximum possible kissing number for $(n - 1)$-dimensional
 Minkowski spheres in $n$-dimensional Euclidean space?
 \end{problem}
 
 \begin{remark}
   The  kissing number of  the Minkowski spheres ${\mathbb S}^{1}_c$, real $c > 1$,  is equal 6 (six).
   This follows from Theorem \ref{tolp}.
 \end{remark}
 
 \begin{conjecture}
  Let $n = 3$.The  kissing number of  the Minkowski spheres ${\mathbb S}^{2}_c$, real $c > 1$, is equal 12 (twelve).\\
Let $n = 4$. The  kissing number of  the Minkowski spheres ${\mathbb S}^{3}_c$, real $c > 1$, is equal 24 (twenty four).
 \end{conjecture}
 
\subsection{Metrizated matroids}. 
\begin{definition}
We will call a matroid {\it metrized} if the metric characteristics of
 its elements are defined, as well as the metric characteristics
 of the mathematical constructions defined on its elements.
\end{definition}

\begin{example}
Let $K$ be a field with a metric on it. Any linear matroid over $K$ is metrizable by extending the field metric to the elements of the matroid.
\end{example}

\begin{proposition}
 \label{pmm}
 Each Minkowski sphere (curve) 
${\mathbb S}^{1}_c$, real $c > 1$, for a given real $c$ parametrizes a continuous family of metrized matroids whose elements are the points of admissible lattices lying on the given Minkowski curve, each having 3 pairs of points lying on the
 curve. Bases of these lattices   can be taken as the bases of matroids. 
 \end{proposition}
 
 Recall some results from \cite{gefms}.
 
  \begin{proposition}
  (Proposition 2 from \cite{gefms})
For Minkowski sphere (curve) 
${\mathbb S}^{1}_c$, real $c > 1$, 
we have next expressions for critical determinants and their lattices:  
\begin{enumerate} 

\item  ${\Delta^{(0)}_p} = \Delta(p, {\sigma_p}) =  \frac{1}{2}{\sigma}_{p},$ 

\item $ {\sigma}_{p} = (2^p - 1)^{1/p},$

\item  ${\Delta^{(1)}_p}  = \Delta(p,1) = 4^{-\frac{1}{p}}\frac{1 +\tau_p }{1 - \tau_p}$,   

\item  $2(1 - \tau_p)^p = 1 + \tau_p^p,  \;  0 \le \tau_p < 1.$  

\end{enumerate}
\end{proposition}
 
  For their critical lattices respectively  $\Lambda_{p}^{(0)},\; \Lambda_{p}^{(1)}$ next conditions satisfy:   $\Lambda_{p}^{(0)}$ and 
 $\Lambda_{p}^{(1)}$  are  two $D_p$-admissible lattices each of which contains
three pairs of points on the boundary of $D_p$  with the
property that 
\begin{itemize}

\item $(1,0) \in \Lambda_{p}^{(0)},$

\item $(-2^{-1/p},2^{-1/p}) \in \Lambda_{p}^{(1)},$

\end{itemize}
 (under these conditions the lattices are
uniquely defined).

Lattice $\Lambda_{p}^{(0)}$ are two-dimensional lattices in ${\mathbb R}^2$ spanned by the vectors
 \begin{itemize}
 \item[]  $\lambda^{(1)} = (1, 0),$
   \item[] $\lambda^{(2)} = (\frac{1}{2}, \frac{1}{2} \sigma_p).$
   \end{itemize}
   For lattices $\Lambda_{p}^{(1)}$, due to the cumbersomeness of the formulas for the coordinates of basis vectors of the
    lattices, as an example, we present only the lattice $\Lambda_{2}^{(1)}$.
    The lattice $\Lambda_{2}^{(1)}$ is a two-dimensional lattice in ${\mathbb R}^2$ spanned by the vectors
 \begin{itemize}
 \item[]  $\lambda^{(1)} = (-2^{-1/2},2^{-1/2}),$
   \item[] $\lambda^{(2)} = (\frac{\sqrt 6 - \sqrt 2}{4}, \frac{\sqrt 6 + \sqrt 2}{4}).$
   \end{itemize}
 
 \begin{lemma}  \cite{gefms}.
Let $(P_x, P_y)$ be a point of the critical lattice $\Lambda$ on Minkowski curve ${\mathbb S}^{1}_{p} = C_{p}$.
Then the point $(u, v)$ that satisfies conditions
\begin{equation}
\label{es}
\left\{
                   \begin{array}{lc}
   |P_x v - P_y u| = \; d(\Lambda),\\
    |u|^p +  |v|^p = \; 1.\\
 \end{array}
                       \right.
\end{equation}
belongs to $\Lambda$ and   lies on $C_{p}$. 
The shell of points of the critical lattice $\Lambda_{p}$ on Minkowski curve $C_{p}$ contains 6 points
Point coordinates $(u, v)$ can be calculated in closed form or with any precision.
\end{lemma}

\begin{example}
\label{ex3}
Each shell of points of the critical lattices $\Lambda_{p}^{(0)}$ on Minkowski curves $C_{p}$ contain 6 points:
\[
  \pm (1, 0), 
  \]
  \[
  \left(\pm \frac{1}{2},  \pm \frac{1}{2} \sigma_p \right).
\]

Respectively the shell of points of the critical lattice $\Lambda_{2}^{(1)}$ on Minkowski curves $C_{p}$ contains 6 points:
\[
  (-2^{-1/2},2^{-1/2}), (2^{-1/2},-2^{-1/2}),  
  \]
  \[
  \pm \left(\frac{\sqrt 6 + \sqrt 2}{4}, \frac{\sqrt 6 - \sqrt 2}{4}\right), 
  \]
  \[
    \pm \left(\frac{\sqrt 6 - \sqrt 2}{4}, \frac{\sqrt 6 + \sqrt 2}{4} \right).
\]

\end{example}

\begin{proposition}
 \label{pmm3}
 Each Minkowski sphere 
${\mathbb S}^{2}_c$, real $c > 1$, for a given real $c$ parametrizes a continuous family of metrized matroids whose elements are the points of admissible lattices lying on the given Minkowski sphere, each having 6 pairs of points lying on the
 sphere. Bases of these lattices   can be taken as the bases of matroids. 
 \end{proposition}

 \subsubsection{Metrizated matroids of  Minkowski spheres}
 For $n = 2$, the metrizated matroids of the  Minkowski spheres ${\mathbb S}^{1}_c$, real $c > 1$, take values in admissible lattices of the spheres ${\mathbb S}^{1}_c$ that contain
three pairs of points on the sphere for a given value of $c$.

\begin{example}
\label{ex4}
  $M({\mathbb S}^{1}_3) = (E,{\mathcal I}), \: \#E = 2, E = \{a, b\}, {\mathcal I} = ab$, where 
  $a = (1,0), \; b =(\frac{1}{2}, \frac{1}{2}{\sigma}_{3} )$. 
  Vectors $a,\: b$ form a basis of the admissible lattice with 
determinant equal to $\frac{1}{2} {\sigma}_{3} = \frac{1}{2} \sqrt[3]{7}$.
\end{example}

 For $n = 3$, the metrizated matroids of the  Minkowski spheres ${\mathbb S}^{2}_c$, real $c > 1$, take values in admissible lattices of the spheres ${\mathbb S}^{2}_c$ that contain
six pairs of points on the sphere for a given value of $c$.

\begin{example}
\label{ex5}
  $M({\mathbb S}^{2}_3) = (E,{\mathcal I}), \: \#E = 3, E = \{a, b,c \}, {\mathcal I} = abc$,  where 
  $a = (1,0,0), \; b =(\frac{1}{2}, \frac{1}{2}{\sigma}_{3}, 0 ), \; c = (0,0,1)$. 
  Vectors $a,\: b, \: c$ form a basis of the admissible lattice with 
determinant equal to $\frac{1}{2} {\sigma}_{3} = \frac{1}{2} \sqrt[3]{7}$.
\end{example}



\begin{thebibliography}{20}

\bibitem{achk}
 Anderson, T., Cook, B., Hughes, K., Kumchev, A. On the ergodic Waring–Goldbach problem, J. Funct. Anal. 282, 39,  Paper No. 109334, 2022.

\bibitem{akc} Arkhipov G.I., Karatsuba A.A., Chubarikov V.N. Theory of multiple trigonometric
 sums. M.: Nauka, 1987.


\bibitem{az}   G. I. Arkhipov, A. N. Zhitkov, On Waring's problem with nonintegral exponent, Math. USSR-Izv., 25:3 , 443–454, (1985).

\bibitem{bb} 
Bannai E.,  Bannai E., A survey on spherical designs and algebraic combinatorics
on spheres,
European Journal of Combinatorics 30  1392-1425 (2009).

\bibitem{bdmp}
Baker M.,  DeMarco L., Preperiodic points and unlikely intersections, Duke Mathematical Journal, 159 (2011), 1-29.

\bibitem{brpt}
 Baker M.,  Rumely R., Potential Theory and Dynamics on the Berkovich Projective Line, Mathematical Surveys and Monographs, 159, American Mathematical Society, 2010.

\bibitem{ber}
Berkovich V. Spectral Theory and Analytic Geometry over Non-Archimedean Fields, American Math Society (1990).

\bibitem{bur}
Bourbaki, Nicolas, Topological Vector Spaces,
 Éléments de mathématique, Springer-Verlag, Berlin New York, 1987.

\bibitem{bmsw}
  Bourgain, J., Mirek, M., Stein, E.M., Wróbel, B.  Dimension-free estimates for discrete HardyLittlewood averaging operators over the cubes in ${\mathbb Z}^d$, Am. J. Math. 141, 857–905, 2019.
  
  \bibitem{bou}
   Bourgain , J. Pointwise ergodic theorems for arithmetic sets, with an appendix by the author, H.
Furstenberg, Y. Katznelson, and D.S. Ornstei, Inst. Hautes Etudes Sci. Publ. Math. 69, 5--45, 1989.

\bibitem{Cassels}
Cassels J. W. S., An Introduction to the Geometry of Numbers, Springer, NY, 1997.

\bibitem{crarot}
H. H. Crapo and G.-C. Rota, On the foundations of combinatorial theory: Combinatorial geometries, preliminary ed., The
M.I.T. Press, Cambridge, Mass.-London, 1970.

\bibitem{Co:MC} H. Cohn, \textit{Minkowski's conjectures on critical lattices
in the metric $\{\vert\xi \vert^p+\vert\eta \vert^p \}^{{1}/{p}},$}
{\it Annals of Math.}, {\bf 51}, (2) (1950), 734--738.

\bibitem{ckmrv}
H. Cohn, A. Kumar, S. D. Miller, D. Radchenko, and M. Viazovska, The sphere
packing problem in dimension 24, Ann. of Math. (2) 185, no. 3, 1017--1033, (2017).

\bibitem{ctf}
Connes, A.,
Trace formula in noncommutative geometry and the zeros of the Riemann zeta function,
Sel. Math. New Ser. 5 (1),  29--106, 1999.

\bibitem{ccrrz}
Connes, A., Consani, C.,
 Riemann–Roch for the ring ${\mathbb Z}$,
 Comptes Rendus. Mathématique.  Vol. 362, p. 229-235, 2024.

\bibitem{ccsp}
Connes, A., Consani, C.,
 Riemann-Roch for $\overline{Spec \; {\mathbb Z}}$,
 Bull. Sci. Math. 187,  1-29, 2023.

\bibitem{ccag}
Connes, A., Consani, C.,
On absolute algebraic geometry the affine case, 
Adv. Math. 390, Article ID 107909, 44 p. 2021.

\bibitem{ccaa}
Connes, A., Consani, C., 
Absolute algebra and Segal’s $\Gamma$-rings,
J. Number Theory, 162, 518-551, 2016.

\bibitem{ccars}
Connes, A., Consani, C.,i
The arithmetic site,
C. R. Math., Ser. I, 352, pp. 971-975, 2014.

\bibitem{ccaar}
Connes, A., Consani, C.,
 From monoids to hyperstructures: in search of an absolute arithmetic, in: Casimir Force, Casimir Operators and the Riemann
Hypothesis, de Gruyter,  pp. 147–198, 2010.

\bibitem{CS}
J.H. Conway, N.J.A. Sloane, 
 Sphere Packings, Lattices and Groups, Second Edition, 
Springer-Verlag, New York Berlin. 1992.

\bibitem{ch}
Cook, B., Hughes, K. Bounds for lacunary maximal functions given by Birch-Magyar averages, Trans. Am. Math. Soc. 374, 3859–3879, 2021.

\bibitem{dei}
Deitmar, A., Schemes over  ${\mathbb F}_1$, in van der Geer, G.; Moonen,
 B.; Schoof, R. (eds.), Number Fields and Function Fields: Two Parallel
 Worlds, Progress in Mathematics, vol. 239, 2005.
 
 \bibitem{del}
  Delsarte, P. An algebraic approach to the association schemes of coding theory, Philips Res. Rep. Suppl. 10, 1973.

\bibitem{etr}
  Ecalle Jean, 
Six Lectures on Transseries, Analysable Functions and the Constructive Proof of Dulac Conjecture,
NATO ASI Bifurcations and Periodic Orbits of Vector Fields (Eds: Dana Schlomiuk), vol. 408, Springer Dordrecht, 75-184.

\bibitem{fe}
Ferraguti, A.
A survey on abelian dynamical Galois groups,
Rend. Semin. Mat., Univ. Politec. Torino 80, 41-54 (2022)

\bibitem{dem}
 DeMarco L.,
 Dynamics of rational maps: A current on the bifurcation locus, Mathematical Research Letters, 8 (2001), 57-66.

\bibitem{de} Deshouillers J.-M. 
Probleme de Waring avec exposants non entiers. Bull. Soc. math.
France,  101, fasc. 3, p. 285--295, 1973. 

\bibitem{dgm}
  Dundas B.,  Goodwillie T.,  McCarthy R., The Local Structure of Algebraic K-Theory, Algebra and
Applications, vol. 18, Springer-Verlag London, Ltd., London, 2013.

\bibitem{du}
Durov N., New approach to Arakelov geometry, arXiv:0704.2030.

\bibitem{edre}
Edelman, P.; Reiner, V. Free hyperplane arrangements between $A_{n-1}$  and $B_n$. (English) 
Math. Z. 215, No. 3, 347-365 (1994).

\bibitem{fkt}
Favre C., Kiwi J., Trucco E., A non-Archimedean Montel's theorem, Compositio Mathematica, 148 (2012), 966-990. 


\bibitem{gl1} Glazunov N.  On packing of Minkowski balls, Comptes rendus de l’Acad´emie bulgare
Sci., Tome 76, No 3, 335-342 (2023).

\bibitem{gl} Glazunov N. On  coverings by Minkowski balls in the plane and a duality,  Comptes rendus de l’Acad´emie bulgare Sci., 
Tome 77, No 6 (2024).

\bibitem{gPacIIl} Glazunov N. 
On packing of Minkowski balls. II,
 arXiv:2302.01644v3 [math.NT] 

\bibitem{GGM:PM} N. Glazunov, A. Golovanov, A. Malyshev,
\textit{Proof of Minkowski's hypothesis about the critical determinant of
$|x|^p+|y|^p<1$ domain}, {\it Research in Number Theory 9}. Notes
of scientific seminars of LOMI. {\bf 151} Leningrad: Nauka. (1986), 40--53.

 \bibitem{Gl4}
N. Glazunov, \textit{On A. V. Malyshev’s approach to Minkowski’s conjecture concerning
the critical determinant of the region $|x|^p +|y|^p < 1$ for $p > 1$}, Chebyshevskii Sb., 
Volume 17, Issue 4  (2016), 185–193.


\bibitem{gl5} Glazunov N. M. Arithmetic statistics, probabilities and Langlands correspondence, Int. Conf.
Analytical and Computational Methods in Probability Theory and its Applications (ACMPT2017), Procedings.  Moscow,   220--225 2017.

\bibitem{gl2} Glazunov N.M. On a mathematical machine oriented to the study of Diophantine 
equations.2.
{\it Computations in algebra, combinatorics and number theory. Kiev: Inst. 
of Cyb. Ac.Ukr.SSR.} (1980).

\bibitem{gl3} Glazunov N.M. On Two Moduli Problems Concerning Number of Points and
Equidistribution over Prime Finite Fields, Proc. of the III Intl. conf. “Discrete models in the theory of control systems” ,  Moscow, Dialogue,  Moscow State University,  23--25 (1998).

\bibitem{gl4} Glazunov N. Methods for Justifying Arithmetic Hypotheses and Computer Algebra.
Programming and computer software, no. 3, 2--8 (2006).

\bibitem{gks} Glazunov N.M., Kaluzhnin L.A., Sushchansky V.I. Programming system for solving
 combinatorial problems of modern algebra, Analytical calculations on
Computers and their applications in theoretical physics: Proc. Int. conf. 
Dubna: JINR, P.23--36 (1980).

\bibitem{gkss} Glazunov N.M.,  Kaluzhnin L.A., Stogniy A.A.,  Sushchansky V.I. Issues of
 development of algebraic research using computers, Cybernetics, no. 2. 1--10, (1983).
 
 \bibitem{gefms} Glazunov N. 
 Extremal functions on moduli spaces and applications,
  arXiv:2411.13671v3 [math.NT]

\bibitem{gefr}
Glazunov N.
Extremal forms and rigidity in arithmetic geometry and in dynamics,
Chebyshevskii Sbornik, 2015, Volume 16, Issue 3, Pages 124--146.

\bibitem{glp}
 Glazjnov N. On Langlands program, global fields and shtukas, Chebyshevskii Sbornik, 2020, Volume 21, Issue 3, Pages 68--83.

\bibitem{glkh}
Glazunov N., Kharchenko V.
Formal and non-Archimedean structures of dynamical systems on manifolds, 
Cybernetics and Systems Analysis, Volume 55, Issue
3 , Page 45-55 (2019).


\bibitem{adr} 
Glazunov N.M.
On algebraic dynamics and resurgence
 on Minkowski moduli spaces, p.48, Book of abstracts,
The 15th Ukraine Algebra Conference, Ivan Franko Nat.Univ.  Lviv, 
July 8-12, 2025.

\bibitem{har} 
  Hartshorne R.,
  Algebraic Geometry,  Springer Science+Business Media, Inc., NY, 1977.

\bibitem{ilms} Iosevich A.,  Langowski B.,  Mirek M., Szarek T.,
     Lattice points problem, equidistribution and ergodic
theorems for certain arithmetic spheres
. (English), Mathematische Annalen  388, No. 2, 2041–2120 (2024).

\bibitem{iss}
   Iosevich, A., Sawyer, E., Seeger, A., On averaging operators associated with convex hypersurfaces of
finite type.  J. Anal. Math. 79, 159–187, 1999.

\bibitem{iro}
Irokawa, R.
Activity measures of dynamical systems over non-Archimedean fields. (English) 
Discrete Contin. Dyn. Syst. 45, No. 2, 361-390 (2025), Zbl 07946836

\bibitem{kra}
  Krasner M., Abstract Galois theory and endotheory I, Acta Sci. Math., 50, 253–286 (1986).

\bibitem{kung}
S. Kung,
 A source book in matroid theory (Boston, Basel, Stuttgart: Birkhäuser) (1986).


\bibitem{lek}
 Lekkerkerker C.G., Geometry of Numbers, NorthHolland, 1969.

\bibitem{kkg}
Krivonos Yu. G., Kharchenko V. P., Glazunov N. M. Differential-algebraic equations and dynamical systems on manifolds. Cybernetics and Systems Analysis, vol. 52, issue 3, 2016, 83-96.


 
 \bibitem{lin}
  Linnik Yu.V.
 Ergodic properties of algebraic fields. Translated from the Russian by M. S. Keane. Ergebnisse der Mathematik und ihrer Grenzgebiete, Band 45 Springer-Verlag New York Inc., NY, 1968.
 
  \bibitem{msw}
 Magyar, Á., Stein, E.M., Wainger, S. Discrete analogues in harmonic analysis: spherical averages. Ann. Math. 155, 189–208, 2002.
 
 \bibitem{ma}
 Manin  Yu., Lectures on zeta functions and motives. Astérisque,
 tome 228, p. 121-163, 1995.

 \bibitem{moow}
Mayhew, Dillon; Oporowski, Bogdan; Oxley, James; Whittle, Geoff
The excluded minors for the class of matroids that are binary or ternary. (English) 
Eur. J. Comb. 32, No. 6, 891-930 (2011). Zbl 1248.05031


\bibitem{mi} W. H. Mills, A prime-representing function, Bull. Amer. Math. Soc. 53(6), 
604--604, 1947.

\bibitem{Mi:DA} Minkowski  H.  {\it Diophantische
Approximationen}, Leipzig: Teubner (1907).

\bibitem{mus}
  Musin O. R.
  The kissing number in four dimensions, Annals of Mathematics, 168, 1–32, (2008).

\bibitem{nak}
  T. Nakasawa, Zur axiomatik der linearen abh¨angigkeit. i, Science Reports of the Tokyo Bunrika Daigaku, Section A 2
(1935), no. 43, 235–255.

\bibitem{neu}
  Neukirch  J. Algebraic number theory, Springer-Verlag, Berlin, Heidelberg, 1999. 

\bibitem{orte} 
 P. Orlik, H. Terao, Arrangements and hypergeometric functions 
(MSJ Memoirs 9, Tokio: Mathematical Society of Japan (2001).

\bibitem{oxl}
Oxley J. On the interplay between graphs and matroids. (English) Zbl 0979.05030
Hirschfeld, J. W. P. (ed.), Surveys in combinatorics, 2001. Proceedings
 of the 18th British combinatorial conference, University of Sussex, 
Brighton, UK, July 1-6, 2001. Cambridge Cambridge University Press. Lond.
 Math. Soc. Lect. Note Ser. 288, 199-239 (2001). Zbl 0979.05030

\bibitem{oxl1}
Oxley J. Matroid theory,  Oxford Science Publications. Oxford Oxford University Press (1992).

\bibitem{pon}
Pontryagin L. S. Ordinary Differential Equations, Elsevier (2014).

\bibitem{rado}
R. Rado, A theorem on independence relations, Quart. J. Math., Oxford Ser., 13,
83--89 (1942).

\bibitem{wiki}
$https://en.wikipedia.org/wiki/Closure_operator$
  
 \bibitem{saffkui}
 Saff E., Kuijlaars A.,
 Distributing many points on a sphere, The Math. Intelligencer, 19, 1, 5--11, 1997.

 \bibitem{seg} B.I. Segal, Waring's theorem for powers with fractional and irrational exponents, Tr. Phys.-math. Institute named after V. A. Steklova, 1934, volume 5, 73–86.

\bibitem{smi}
Smirnov A., Hurwitz inequalities for number fields, Algebra i 
Analiz (in Russian), 4 (2), 186--209, 1992.

\bibitem{sou}
Soulé  C., Les variétés sur le corps à un élément, Mosc. Math. J. 4 (1),  21, 2004.

\bibitem{stei}
N. Steenrod, S. Eilenberg,
 Foundations of algebraic topology,  Princeton university press, Princeton, 1952.

\bibitem{sta}
R. Stanley, Modular elements of geometric lattices, Algebra Univers.
1, 214-217 (1971).

\bibitem{sta1} 
R/ Stanley, Supersolvable Lattices, Algebra Univers. 2, 197-217 (1972)

 
 \bibitem{tit}
 Tits J., Sur les analogues algébriques des groupes semi-simples complexes,
 Colloque d'algèbre supérieure, tenu à Bruxelles du 19 au 22 décembre 1956, 
Centre Belge de Recherches Mathématiques Établisements Ceuterick, Louvain, 
Paris: Librairie Gauthier-Villars, pp. 261–289, 1957.

 \bibitem{tva}
Toën B., 
 Vaqui\`e  M., Au dessous de $Spec \; {\mathbb Z}$, J. K-Theory 3 (3),  437–500, 2009.

 \bibitem{tru}
Truemper, K., Matroid Decomposition (1992), Academic Press Inc.: Academic Press Inc. Boston, MA 
 
\bibitem{tut}
Tutte, W. T. Graph Theory, Cambridge University Press, (2001)., 

\bibitem{var}
A. Varchenko, Proc. Int. Congr. Math., Kyoto/Japan 1990, Vol. I, 281-300 
(1991).

\bibitem{vi} I. M. Vinogradov. Method of trigonometric sums in number theory. M.: Nauka,
1980.

\bibitem{vi1} I.M. Vinogradov, On Waring’s theorem, Proceedings of the USSR Academy of
 Sciences. VII series. Department of Physical and Mathematical Sciences, 
1928, No. 4, 393–400.

\bibitem{vi2} I. M. Vinogradov, Representation of an odd number as a sum of three primes,
 Comptes Rendues (Doklady) de l’Academy des Sciences de l’USSR 15 (1937), 
191–294.

\bibitem{via} 
 Viazovska M.S. The sphere packing problem in
dimension 8, Ann. of Math. (2)
185, no. 3, 991–1015, (2017).

\bibitem{wie}
Wielandt H. Permutation representation, J. of Math., 13, no.1, 91-94.   

\bibitem{whi}
 Whitney H. On the abstract properties of linear dependence, Amer. J.
Math. 57, 509–533 (1935).

\bibitem{zas}
Zaslavsky T.  Biased graphs. I. Bias, balance, and gains, JOURNAL OF COMBINATORIAL THEORY, Series B 47, 32-52 (1989). 

\bibitem{zas1}
  Zaslavsky  T. Supersolvable frame-matroid and graphic-lift lattices (English),
Eur. J. Comb. 22, No. 1, 119-133 (2001). Zbl 0966.05013

\end{thebibliography}
\end{document}